\documentclass[bj,authoryear]{imsart}

\RequirePackage{amsmath,amsthm,amssymb,enumerate,amsfonts,amssymb,mathtools,subcaption, bigints, comment,bbm}
\RequirePackage{graphicx}

\startlocaldefs

\numberwithin{equation}{section}
\theoremstyle{plain}

\newtheorem{theorem}{Theorem}[section]
\newtheorem{lemma}[theorem]{Lemma}
\theoremstyle{remark}
\newtheorem{definition}[theorem]{Definition}

\newtheorem{assumption}{Assumption}
\newtheorem*{remark}{Remark}

\newcommand{\uA}       {\mbox{\boldmath$A$}} 
 
\newcommand{\uB}       {\mbox{\boldmath$B$}}

\newcommand{\uD}       {\mbox{\boldmath$D$}}

\newcommand{\uI}       {\mbox{\boldmath$I$}}

\newcommand{\uM}       {\mbox{\boldmath$M$}}

\newcommand{\uO}       {\mbox{\boldmath$O$}}

\newcommand{\uS}       {\mbox{\boldmath$S$}} 

\newcommand{\uT}       {\mbox{\boldmath$T$}}

\newcommand{\uW}       {\mbox{\boldmath$W$}}

\newcommand{\uX}       {\mbox{\boldmath$X$}}

\newcommand{\uY}       {\mbox{\boldmath$Y$}}
\newcommand{\uYn}       {\mbox{\boldmath$Y_n$}}

\newcommand{\uiota}             {\mbox{\boldmath$\uiota$}}

\newcommand{\uDelta}            {\mbox{\boldmath$\Delta$}}

\newcommand{\uSigma}            {\mbox{\boldmath$\Sigma$}}

\newcommand{\uPsi}              {\mbox{\boldmath$\Psi$}}
\newcommand{\uOmega}            {\mbox{\boldmath$\Omega$}}

\newcommand{\uone}               {\mbox{\boldmath$1$}}
\DeclareMathOperator{\tr}{tr}

\newcommand{\norm}[1]{\left\lVert#1\right\rVert_2}
\newcommand{\normr}[1]{\left\lVert#1\right\rVert}
\newcommand{\snorm}[1]{\left\lVert#1\right\rVert_{\psi_2}}

\newcommand{\fnorm}[1]{\left\lVert#1\right\rVert_{F}}

\endlocaldefs

\begin{document}

\begin{frontmatter}
\title{High-dimensional Bernstein Von-Mises theorems for covariance and precision matrices}
\runtitle{BvM theorems for covariance and precision matrices}


\begin{aug}
\author[A]{\fnms{Partha}~\snm{ Sarkar$^\ast$}\ead[label=e1]{psarkar@fsu.edu}}
\author[B]{\fnms{Kshitij}~\snm{Khare}\ead[label=e2]{kdkhare@stat.ufl.edu}}
\author[B]{\fnms{Malay}~\snm{Ghosh}\ead[label=e3]{ghoshm@ufl.edu}}
\author[C]{\fnms{Matt}~\snm{P. Wand}\ead[label=e4]{Matt.Wand@uts.edu.au}}
\address[A]{Department of Statistics,
Florida State University\printead[presep={,\ }]{e1}}
\address[B]{Department of Statistics,
University of Florida\printead[presep={,\ }]{e2,e3}}
\address[C]{School of Mathematical and Physical Sciences,
University of Technology Sydney,
\printead{e4}}
\end{aug}

\begin{abstract}
This paper aims to examine the characteristics of the posterior distribution of covariance/precision matrices in a ``large $p$, large $n$" scenario, where $p$ represents the number of variables and $n$ is the sample size. Our analysis focuses on establishing asymptotic normality of the posterior distribution of entire covariance/precision matrices under specific growth restrictions on $p_n$ and other mild assumptions. In particular, the limiting distribution is a symmetric matrix variate normal distribution whose parameters depend on the maximum likelihood estimate. Our results hold for a wide class of prior distributions which includes standard choices used by practitioners. Next, we consider Gaussian graphical models that induce precision matrix sparsity. The posterior contraction rates and asymptotic normality of the corresponding posterior distribution are established under mild assumptions on the prior and true data-generating mechanism.
\end{abstract}


\end{frontmatter}


\section{Introduction}\label{sec1}

The advent and proliferation of high-dimensional data and associated Bayesian statistical methods in recent years have generated significant interest in establishing high-dimensional asymptotic guarantees for such methods. The Bernstein von-Mises (BvM) theorem ( \cite{bernstein,cam2000asymptotics, vaart, vonMises}) is a key result that can justify Bayesian methods from a frequentist point of view. The BvM approach assumes a frequentist data-generating model and defines criteria for the prior that result in the posterior becoming asymptotically Gaussian as the number of observations $n$ increases. The primary use of the BvM method is to justify the construction of Bayesian credible sets as a Bayesian counterpart of the frequentist confidence region. It is useful in cases where uncertainty quantification through frequentist methods is not feasible due to the presence of unknown parameters in the asymptotic distribution, making it challenging to construct frequentist confidence regions directly.

Although there is extensive literature establishing Bernstein von Mises theorems in settings where the number of parameters $p$ stays fixed as $n$ increases, analogous results for high-dimensional settings where $p = p_n$ can grow with sample size $n$ are comparatively sparse. 
In the context of linear models, BvM results were established by \cite{Dominique, castillo1, ghosal1}, while \cite{Boucheron, clarke, ghosal3, ghosal2} studied it for high-dimensional exponential models, subject to certain conditions on the growth rate of the dimension. \cite{spokoiny} explored similar ideas in a wider ``general likelihood setup". \cite{panov} explored BvM results in a semiparametric framework with finite sample bounds for distance from normality since modern statisticians are increasingly focused on models with limited sample sizes. See also \cite{kelijn, castillo2, rivoirard, spokoiny1,spokoiny2, katsevich1, katsevich2} for additional results in this context.

Our focus in this paper is Bayesian methods for high-dimensional covariance estimation. In particular, suppose we have $n$ independent and identically distributed samples $\uY^{n}=(Y_1,\dots,Y_n)$ drawn from a $p$-variate normal distribution with covariance matrix $\uSigma$. We first consider the ``unstructured" estimation of $\uSigma$, i.e., no dimension-reducing structure, such as sparsity or low-rank, is imposed on $\uSigma$. In this setting, \cite{chaogao} studied BvM results for one-dimensional functionals of the covariance matrix such as matrix entries and eigenvalues, in a high-dimensional setting. \cite{silin} derived finite sample bounds for the total variation distance between the posterior distributions of $\uSigma$ obtained by employing an Inverse-Wishart (IW) prior and a flat prior. Moreover, he investigated Bernstein-von Mises theorems for one-dimensional functionals and spectral projectors of the covariance matrix.

However, when it comes to simultaneously inferring various functionals of the covariance matrix (such as multiple entries of the covariance matrix, its inverse, and multiple eigenvalues) the above results are not applicable even with a very basic conjugate family of Inverse-Wishart (IW) priors. Although a Bonferroni inequality-based approach could potentially be utilized, it often results in inefficient and loose bounds, particularly in high-dimensional settings. The key goal of this paper is to provide a high dimensional Bernstein-von Mises theorem for the entire covariance matrix $\uSigma$ (or the precision matrix $\uOmega$) for a general enough class of priors. We show that as long as the prior distribution satisfies the flatness condition around sample covariance matrix $\uS$ (see equation \ref{flat}) the total variation norm between the posterior distribution of $\sqrt{n}(\uSigma-\uS)$ (or $\sqrt{n}(\uOmega-\uS^{-1})$) and a suitable mean zero symmetric matrix variate normal distribution tends to zero under standard regularity assumptions (Theorem \ref{th_sigma} and \ref{th_omega}). We show that a large collection of the prior distributions for $\uSigma$ (or $\uOmega$) satisfy this flatness condition around $\uS$ (Lemma \ref{postcont1}, \ref{postcont2} and \ref{postcont3}). This includes standard conjugate IW prior and several scale mixtures of IW prior proposed in \cite{gelman2006prior, gelman_book,gelman2006data, huang, mulder,zava}. These mixtures have been shown to offer more effective noninformative choices. In fact, we are able to show that the flatness condition around $\uS$ is satisfied by a significantly generalized version of the mixture priors proposed in the above literature.

Establishing BvM results for the entire covariance matrix poses a significant challenge, especially in high-dimensional settings. The primary issue arises from the fact that an unrestricted $(p\times p)$ covariance/precision matrix involves a large number of free parameters, which is $O(p^2)$. Consequently, as the dimension increases and $p_n$ grows with $n$, the number of parameters escalates rapidly. Furthermore, as discussed in \cite{ghosal1}, when $p_n$ grows with the sample size $n$, there can be a tail region where the posterior probability is significant, even if the likelihood is small in that region. With these challenges, we establish BvM results for the entire covariance matrix $\uSigma$ where $p(=p_n)$ can increase with $n$ but is subject to the condition that $p_n^{5}=o(n)$ (see Theorems \ref{th_sigma} 
and \ref{th_omega}). This seemingly stringent requirement is not due to any imprecise bounds in the proof and is somewhat expected given related results under simpler settings in the literature. \cite{silin} requires the same condition to establish the asymptotic equivalence (in TV norm) of posterior distributions using IW prior and a flat prior. In a simpler context of BvM results for high-dimensional regression condition $p_n^4(\log(p_n))=o(n)$ is required in \cite{ghosal1}. To establish BvM results for several \textit{one-dimensional} functionals of $\uSigma$, the authors in \cite{chaogao} need the condition $p_n^4=o(n)$. 

Recall that the above discussion focuses on a setting where no structure is imposed on the covariance matrix to reduce its dimensionality. A standard and popular approach for parameter reduction in high-dimensional covariance estimation settings is to impose sparsity in the precision matrix. These models are referred to as Gaussian graphical models or concentration graphical models (see \cite{lauritzen}). A specific sparsity pattern in $\uOmega$ can be conveniently represented by a graph $G$ involving the set of $p$ variables. The $G$-Wishart distribution, as introduced by \cite{roverato}, offers a conjugate family of priors for the concentration graphical model corresponding to a given graph $G$. Decomposable graphs, which have received considerable attention in Bayesian literature on concentration graph models (see \cite{dawid,letac, rajaratnam, roverato}), form a notable subfamily within this framework. High-dimensional posterior contraction rates for the precision matrix in these models have been established in \cite {banerjee1,xiang} (underlying graph known) and \cite{liu2019empirical, Lee} (underlying graph unknown). These contraction rates play a crucial role in establishing BvM results, and hence it is important to ensure their sharpness/optimality. However, either these posterior contraction rates are not close to the optimal frequentist convergence rates (see \cite{rothman}) established in the literature, or they require stringent conditions that render them inapplicable in high-dimensional settings, even when the underlying graph is known.

 We address this issue by establishing Frobenius norm posterior contraction rates for the precision matrix (under a decomposable concentration graphical model) which match the optimal frequentist convergence rates in \cite{rothman} for both cases when the underlying true graph is known (refer to Theorem \ref{postthm2}) or unknown (refer to Theorem \ref{postthm4} and \ref{postthm5}). Additionally, we establish posterior contraction rates under the spectral norm (see Theorem \ref{postthm1} and Theorem \ref{postthm3}) which significantly improve previous rates in \cite{xiang,banerjee1}. Leveraging these posterior contraction rates, we derive, under mild regularity conditions, a BvM result for the precision matrix when the imposed sparsity pattern corresponds to a decomposable graph for both cases when the underlying true graph is known (refer to Theorem \ref{th_omega_sp}) or unknown (refer to Theorem \ref{th_omega_sp1}). If the maximum vertex degree of this graph is assumed to be bounded (e.g. see \cite{banerjee1}), then the condition $p_n^5=o(n)$ that we needed for the unstructured setting is significantly weakened to  $p_n^2(\log(p_n))^3=o(n)$.

Section \ref{sec11} of the paper aims to demonstrate that there is nothing special about using total variation norms for BvM results. Other distance measures, such as the Bhattacharyya-Hellinger distance (\cite{bhattacharyya1946measure, Hellinger}) or R\'enyi's $\alpha$-divergence (\cite{renyi}), can also be employed to draw similar conclusions.

The remainder of the paper is organized as follows. After introducing the basic notation in the next subsection, the fundamental definitions and preliminaries are presented in Section~\ref{sec2}. Section~\ref{sec3} discusses various prior distributions for dense covariance or precision matrices, and the BvM results for this unstructured dense setting are given in Section~\ref{sec4}. Preliminaries related to concentration graphical models appear in Section~\ref{sec7}. Sparsity-based models for the precision matrix and the corresponding prior distributions are formulated in Section~\ref{sec8}. The BvM and posterior consistency results for the case of a known underlying graph structure are provided in Section~\ref{sec9}, while analogous results for the unknown-graph setting are presented in Section~\ref{sec10}. Section~\ref{sec11} addresses the equivalence of various matrix norms in the context of convergence. Proofs of selected theorems and technical lemmas are given in the supplementary document \cite{supp1}. Finally, we conclude the paper with a summary of our findings and closing remarks in Section~\ref{sec11}.

\subsection{Notation}\label{1.1}
Let us introduce some notation and definitions. For positive sequences $a_n$ and $b_n$, we denote $a_n = O(b_n)$ if there exists a constant $C$ such that $a_n \leq C b_n$ for all $n \in \mathbb{N}$, and $a_n = \Omega(b_n)$ if there exists a constant $C$ such that $a_n \geq C b_n$ for all $n \in \mathbb{N}$. We use $a_n = o(b_n)$ to denote the limit $\lim_{n \to \infty} \frac{a_n}{b_n} = 0$. Also, $a_n \sim b_n$ means that $\frac{a_n}{b_n} \to 1$ as $n \to \infty$. Given a metric space $(X, d),\;\mathcal{N}(X,\epsilon)$ represents the $\epsilon$-covering number, which is the minimum number of balls of radius $\epsilon$ needed to cover $X$.

In the following notation, $\uI_p$ represents the identity matrix of order $p$, and $\uO$ represents a matrix of size $a \times b$ with all zero entries. If $\uA$ is a symmetric square matrix, then $\lambda_{\min}(\uA)$ and $\lambda_{\max}(\uA)$ denote the smallest and largest eigenvalues of $\uA$, respectively. The tensor or Kronecker product between two matrices $\uA$ and $\uB$ is denoted by $\uA\otimes\uB$. Consider the set $M_p$, which comprises all symmetric matrices of size $p \times p$, and a subclass of $M_p$, $\mathbb{P}_p^+$, representing the collection of symmetric positive definite matrices of size $p \times p$.

The unit Euclidean sphere in $\mathbb{R}^p$ is denoted by $\mathcal{S}^{p-1}$. For a vector $x \in \mathbb{R}^p$, we denote its $r$-th norm by $\normr{x}_r = \left(\sum_{j=1}^p \lvert x_j\rvert^r\right)^{1/r}$. $\norm{x}$ denotes the Euclidean norm. For a $p\times p$ matrix $\uA=(A_{ij})_{1\leq i,j\leq p}$, we denote
$$\normr{\uA}_{\max}=\underset{1\leq i,j\leq p}{\max}\lvert A_{ij}\rvert, $$
$$\normr{\uA}_{r,s}=\sup\left\{\normr{\uA x}_{s}:\normr{x}_{r}=1\right\},$$
where $1\leq r,s\leq \infty$. In particular, we have $\normr{\uA}_{(\infty,\infty)}=\underset{i}{\max}\sum_j\lvert A_{ij}\rvert$, and spectral norm of a matrix is defined as $$\norm{\uA}:=\sup_{u\in\mathcal{S}^{n1}}\norm{\uA u}(=\normr{\uA}_{2,2}).$$ 
We define the vectorization of $\uA$ as $\mbox{vec}(\uA) = (A_{11},\dots,A_{p1},A_{12},\dots,A_{pp})^T$. If $\uA$ is a symmetric matrix, there will be repeated elements in $\mbox{vec}(\uA)$. For a $p \times p$ symmetric matrix $\uA$, $\mbox{vech}(\uA)$ is a column vector of dimension $\frac{1}{2}p(p+1)$ formed by taking the elements below and including the diagonal, column-wise. In other words, $\mbox{vech}(\uA) = (A_{11},A_{21},\cdots,A_{p1},A_{22}\cdots,A_{p2},\cdots,A_{pp})^T$. For a symmetric matrix $\uA$, we can establish the connection between $\mbox{vec}(\uA)$ and $\mbox{vech}(\uA)$ using an elimination matrix $\uB_p^T$, expressed as $\mbox{vech}(\uA) = \uB_p^T \mbox{vec}(\uA)$. While it's important to note that this elimination matrix may lack uniqueness, we can construct a $\frac{1}{2}p(p+1) \times p^2$ elimination matrix $\uB_p^T$ in the following systematic manner, as described in \cite{magnus}
$$\uB_p^T=\underset{1\leq j \leq i \leq p}{\sum}\left(u_{ij}\otimes e_j^T\otimes e_i^T\right),$$
where $e_{i}$ is a unit vector whose $i$-th element is one and zeros elsewhere and $u_{ij}$  is a unit vector of order $\frac{1}{2}p(p+1)$ having the value 
$1$ in the position $(j-1)p+i-{\frac {1}{2}}j(j-1)$ and $0$ elsewhere.

Let $f$ and $g$ be two densities, each continuous with respect to some $\sigma$-finite measure $\mu$. Also, let $P(A) =\int_A f d\mu$ and $Q(A) =\int_A g d\mu$. Then Total Variation (TV) norm between two distributions $P$ and $Q$ (or two densities $f$ and $g$) is defined as
$$TV(f,g)=\underset{A}{\sup}\;\lvert P(A)-Q(A)\rvert=\frac{1}{2}\int\lvert f-g\rvert\;d\mu.$$

\section{Preliminaries and model formulation in the dense setting}\label{sec2}

We consider an independent and identically distributed
sample of size $n$, $\uY^n=(Y_1,\dots,Y_n)$, drawn from the
$N_p(0,\;\uSigma=\uOmega^{-1})$ distribution. In the case of estimating $\uSigma$, the moment-based or maximum likelihood estimation is represented as $\uS=\frac{1}{n}\sum_{i=1}^{n}Y_{i}Y_{i}^T$, whereas, for $\uOmega$, it is $\uS^{-1}$. Within this Gaussian framework, we can express the log-likelihood function of $\uSigma$, denoted as $l_{1n}(\uSigma)$, as follows
\begin{align}\label{likelihood1}
    l_{1n}(\uSigma)=-\frac{np}{2}\log(2\pi)-\frac{n}{2}\log(\det(\uSigma))-\frac{n}{2}\tr(\uSigma^{-1}\uS)
\end{align}
Similarly, we can write the log-likelihood function of $\uOmega$, denoted as $l_{2n}(\uOmega)$, as
\begin{align}\label{likelihood2}
    l_{2n}(\uOmega)=-\frac{np}{2}\log(2\pi)+\frac{n}{2}\log(\det(\uOmega))-\frac{n}{2}\tr(\uOmega \uS)
\end{align}
In a Bayesian framework, a prior $\Pi_{1n}(\cdot)$ is assigned to the covariance matrix $\uSigma$. The induced prior on the precision matrix $\uOmega$ is denoted as $\Pi_{2n}(\cdot)$. Let $\pi_{1n}(\cdot)$ and $\pi_{2n}(\cdot)$ represent the corresponding prior densities. We will consider an asymptotic framework where $p$ will be allowed to grow with the sample size $n$. This is why the dependence of the priors on $n$ is highlighted in the above notation. For the sake of notational simplicity, we will sometimes refer to $\pi_{1n}(\cdot)$ and $\pi_{2n}(\cdot)$ as $\pi_{1}(\cdot)$ and $\pi_{2}(\cdot)$, respectively.

Now, after centering $\uSigma$ by $\uT_1=\sqrt{n}(\uSigma-\uS)$ or $\uOmega$ by $\uT_2=\sqrt{n}(\uOmega-\uS^{-1})$, we define the following functions
\begin{align}
&M_{1n}(\uT_1)=\exp{\left(l_{1n}\left(\uS+\frac{\uT_1}{\sqrt{n}}\right)-l_{1n}\left(\uS\right)\right)},\;\text{and}\label{m1}\\
&M_{2n}\left(\uT_2\right)=\exp{\left(l_{2n}\left(\uS^{-1}+\frac{\uT_2}{\sqrt{n}}\right)-l_{2n}\left(\uS^{-1}\right)\right)}\label{m2},
\end{align}
\noindent where $\uT_1 \in B_{1n}$ and $\uT_2 \in B_{2n}$, where $B_{1n}=\{\uT_1:\uS+\dfrac{\uT_1}{\sqrt{n}}\in\mathbb{P}_p^+\}$ and $B_{2n}=\{\uT_2:\uS^{-1}+\dfrac{\uT_2}{\sqrt{n}}\in\mathbb{P}_p^+\}$. If $\uT_1$ or $\uT_2$ falls outside $B_{1n}$ or $B_{2n}$, we set $M_{1n}(\uT_1)$ or $M_{2n}(\uT_2)$ equal to zero. Suppose, posterior distributions of $\uT_1$ and $\uT_2$ are given by $\Pi_{1n}(\cdot\mid\uY^n)$ and $\Pi_{2n}(\cdot\mid\uY^n)$ respectively. Analogously, let $\pi_{1n}(\cdot\mid\uY^n)$ and $\pi_{2n}(\cdot\mid\uY^n)$ be the corresponding posterior densities. Then it is not difficult to check
\begingroup
\allowdisplaybreaks
\begin{align}
    &\pi_{1n}(\uT_1\mid\uY^n)=\frac{M_{1n}(\uT_1)\pi_{1}\left(\uS+\dfrac{\uT_1}{\sqrt{n}}\right)}{\bigintsss_{B_{1n}} M_{1n}(\uW)\pi_{1}\left(\uS+\dfrac{\uW}{\sqrt{n}}\right)d\uW}\;, \text{and}\label{post1}\\
    &\pi_{2n}(\uT_2\mid\uY^n)=\frac{M_{2n}(\uT_2)\pi_{2}\left(\uS^{-1}+\dfrac{\uT_2}{\sqrt{n}}\right)}{\bigintsss_{B_{2n}} M_{2n}(\uW)\pi_{2}\left(\uS^{-1}+\dfrac{\uW}{\sqrt{n}}\right)d\uW}.\label{post2}
\end{align}
\endgroup
Before proceeding with further discussion, let us revisit two important definitions from the literature.
\begin{definition}{\textbf{(Symmetric Matrix-normal Distribution)}}\label{def1}
Let $\uX$ be a $p\times p$ symmetric random matrix, $\uM$ is a $p\times p$ deterministic symmetric matrix and say $\uPsi_1$ and $\uPsi_2$ be constant $p\times p$ positive definite symmetric matrices such that $\uPsi_1\uPsi_2=\uPsi_2\uPsi_1$. Then $\uX (=\uX^T)$ is said to have a symmetric matrix-normal distribution, denoted by $\mathcal{SMN}_{p \times p}(\uM,\;\uB_p^T(\uPsi_1\otimes\uPsi_2)\uB_p)$, if and only if $\mbox{vech}(\uX)\sim\mathcal{N}_{p(p+1)/2}(\mbox{vech}(\uM),\; \uB_p^T(\uPsi_1\otimes\uPsi_2)\uB_p)$. The probability density function $f(\cdot)$ of a  $\mathcal{SMN}_{p \times p}(\uM,\;\uB_p^T(\uPsi_1\otimes\uPsi_2)\uB_p)$ can be expressed as follows
$$f(X)=\frac{\exp\{-\tr(\uPsi_1^{-1}(X-\uM)\uPsi_2^{-1}(X-\uM))/2\}}{(2\pi)^{p(p+1)/4}\det(\uB_p^T(\uPsi_1\otimes\uPsi_2)\uB_p)^{1/2}},\;X \in M_p.$$
\end{definition}
\noindent 

\noindent If $\uX\sim\mathcal{SMN}_{p \times p}(\uM,\;\uB_p^T(\uPsi_1\otimes\uPsi_2)\uB_p)$ and if $\uA$ is a $q \times p$ matrix of rank $q(\leq p)$ then $\uA\uX\uA^T\sim\mathcal{SMN}_{q \times q}(\uA\uM\uA^T,\;\uB_{q}^T((\uA\uPsi_1\uA^T)\otimes(\uA\uPsi_2\uA^T))\uB_{q})$. For more properties of a symmetric matrix variate normal distribution see \cite{gupta}.
\begin{definition}{\textbf{(Sub-Gaussian Random Variable)}}\label{def2} A mean zero random variable X that satisfies $E[\exp(tX)]\leq \exp(t^2k_1^2)$ for all $t \in \mathbb{R}$ and a constant $k_1$ is called a sub-Gaussian random variable. 
\end{definition}
\noindent
If X is a sub-Gaussian random variable then it satisfies $\left( E(\lvert X\rvert^q)\right)^{1/q}\leq k_2\sqrt{q}$ for a constant $k_2$ and one can define its sub-Gaussian norm as $\snorm{X}\coloneqq \sup_{q\geq1}\;q^{-1/2}\left( E(\lvert X \rvert^q)\right)^{1/q}$. A mean zero random vector $X\in \mathbb{R}^p $ is said to be sub-Gaussian if, for any $u\in \mathcal{S}^{p-1}$, the random variable $u^{T}X$ is sub-Gaussian. The sub-Gaussian norm of a random vector X is defined as,
\begin{align}
    \snorm{X}\coloneqq \underset{u\in \mathcal{S}^{p-1}}{\sup}\;\snorm{u^{T}X}.\notag
\end{align}

\noindent
See \cite{vershynin} for more details.

We now specify the {\it true data-generating mechanism}. As mentioned earlier, we denote the dimensionality of the responses as $p_n$ to highlight the fact that the number of responses, denoted by $p$, can grow with the sample size $n$, making our results applicable to high-dimensional scenarios. We assume that the observations $Y_1,\dots, Y_n$ are independently and identically distributed from a sub-Gaussian random variable with zero mean, where the variance of $Y_1$ is denoted as $\uSigma_{0n}$ (or $\uOmega_{0n}^{-1}$). Thus, the sequence of true covariance (or precision) matrices is represented as $\{\uSigma_{0n}\}_{n\geq1}$ (or $\{\uOmega_{0n}\}_{n\geq1}$). For convenience, we denote $\uSigma_{0n}^{p_n\times p_n}$ as $\uSigma_{0}$ and $\uOmega_{0n}^{p_n\times p_n}$ as $\uOmega_{0}$, specifically highlighting that $\uSigma_{0}$ or $\uOmega_{0}$ depends on $p_n$ (and therefore on $n$). Let $\mathbb{P}_{0n}$ denote the probability measure underlying the true model described above. To simplify notation, we will use $\mathbb{P}_{0}$ instead of $\mathbb{P}_{0n}$. With all this notion in hand, we will define a notion of \textit{posterior consistency} for $\uSigma$ as follows.
\begingroup
\allowdisplaybreaks
\begin{definition}\label{def3}
The sequence of marginal posterior distributions of $\uSigma$ given by $\{\Pi_n(\uSigma\mid\uYn)\}_{n\geq1}$ is said to be consistent at $\uSigma_{0}$, if for every $\delta>0$, $$\Pi_n(\norm{\uSigma-\uSigma_{0}} > \delta \mid \uYn) \overset{P}{\to} 0$$ as $n \to \infty$, under $\mathbb{P}_{0}$. Additionally, if $\Pi_n(\norm{\uSigma-\uSigma_{0}} > \varepsilon_n \mid \uYn) \overset{P}{\to} 0$ as $n \rightarrow \infty$, under $\mathbb{P}_{0}$ for some sequence $\varepsilon_n \rightarrow 0$. Then we refer to $\varepsilon_n$ as the contraction rate of $\{\Pi_n(\uSigma\mid\uYn)\}_{n\geq1}$ around $\uSigma_0$. 
\end{definition}
\endgroup

\noindent Likewise, posterior consistency and contraction rate for $\uOmega$ can also be defined. Also, the posterior consistency and the contraction rate can alternatively be expressed using the Frobenius norm.  The relationship between the concept of posterior consistency and the BvM theorem may not be immediately apparent at this stage. However, in various contexts, posterior consistency is often a crucial requirement for proving BvM results (refer to \cite{chaogao, ghosal1} for specific instances). Further discussion can be found in Section \ref{sec4}. 

With the notion of posterior consistency in hand, similar to \cite{silin} we define the concept of \textit{flatness of a prior around the sample covariance matrix \uS}. For the sequence of priors $\{\Pi_{1n}(\cdot)\}_{n\geq1}$, let us define $\rho^{\pi_{1}}(\varepsilon_n)$ as follows

\begin{align}\label{flat}
\rho^{\pi_{1}}(\varepsilon_n):= \underset{\mathbf{T_1}\in D(\varepsilon_n)}{\sup}\;\left\lvert\frac{\pi_{1}(\uS+\dfrac{\uT_1}{\sqrt{n}})}{\pi_{1}(\uS)}-1\right\rvert,
\end{align}
where $D(\varepsilon_n)=\{\uT_1\in B_{1n}\mid\norm{\uT_1} \leq \sqrt{n}\varepsilon_n \}$. 
Note that, $\uT_1\in D(\varepsilon_n)$ if and only if $\norm{\uSigma-\uS} \leq \varepsilon_n$, where $\varepsilon_n$ is the posterior contraction rate.
A prior distribution with density $\pi_{1}(\cdot)$ is considered flat around $\uS$ if $\rho^{\pi_{1}}(\varepsilon_n)$ tends to $0$ in probability as $n$ approaches infinity. Similarly one can define the flatness of prior around the inverse of sample covariance matrix $\uS^{-1}$ for the sequence of induced prior distribution( $\{\Pi_{2n}(\cdot)\}_{n\geq1}$) for the precision matrix using posterior contraction rates of $\uOmega$. Note that, Definition \ref{def3} considers contraction around the true value $\uSigma_{0}$, hence this notion of flatness will be useful only when $\uS$ also contracts around $\uSigma_{0}$ at the same rate $\varepsilon_n$. Fortunately, this holds for all classes of prior distribution that will be considered in the next section (see the discussion after Assumption \ref{as3}). 
A similar type of condition has also been imposed by \cite{ghosal1, Yano} in the BvM literature. However, their formulation differs slightly, though the underlying essence is the same as ours: they assume that the joint prior distribution is locally log-Lipschitz. In the next lemma, we show that if any joint prior distribution is locally log-Lipschitz with respect to the spectral norm, then the flatness condition is automatically satisfied. So, in that sense, our condition is weaker. The proof of this lemma is provided in the supplementary material (\cite{supp1}). 

\begin{lemma}[Local log-Lipschitz prior implies flatness]
\label{lem:lip_implies_flat}
Suppose that, for each $n$, the prior density $\pi_1(\cdot)$ satisfies the following local log-Lipschitz condition around the sample covariance matrix $S$: there exist deterministic sequences $r_n > 0$ and $L_n > 0$ such that, with probability tending to one under $\mathbb{P}_0$,
\begin{equation}\label{eq:local_log_lip}
  \bigl|\log \pi_1(\Sigma) - \log \pi_1(S)\bigr|
  \;\le\;
  L_n \,\|\Sigma - S\|_2
  \quad
  \text{for all $\Sigma\in \mathbb{P}_p^+$ with }
  \|\Sigma - S\|_2 \le r_n.
\end{equation}

\noindent Assume further that the posterior contraction rate $\epsilon_n$ and the neighborhood radius
$r_n$ satisfy $ \varepsilon_n \;\le\; r_n$ and $L_n \epsilon_n \;\longrightarrow\; 0,$ as $n \to \infty$. Then the flatness condition around $S$ holds, i.e.
\begin{equation}\label{eq:flatness_conclusion}
  \rho_{\pi_1}(\epsilon_n)
  \;\xrightarrow{\mathbb{P}_0}\; 0
  \quad\text{as } n \to \infty.
\end{equation}
\end{lemma}

\noindent In our analysis we have observed that, for the unstructured covariance (or precision matrix) case with an Inverse–Wishart (or Wishart) prior, a bound of the form $L_n = O(p_n^{2})$ (see Lemmas~\ref{postcont1}, \ref{postcont2}, \ref{postcont3}) is sufficient to ensure prior flatness around $S$. In contrast, for the sparse precision matrix setting under a $G$–Wishart prior, a substantially weaker requirement of $L_n = O(p_n)$ (see Lemmas~\ref{gw} and \ref{gw1}) guarantees the same property. One might wonder why a weaker condition such as $L_n = O(p_n)$ was adequate in \cite{Yano, ghosal1} even in the unstructured parameter case. The key point is that their setting assumes either i.i.d.\ or multivariate normal prior distributions on the entries, whereas in our framework we work with significantly more intricate Wishart-type priors, where the determinant term introduces a nontrivial layer of dependence. This added structural complexity necessitates a stronger growth condition on $L_n$ to verify the required flatness property. The key idea behind such \textit{log-Lipschitz} or \textit{flatness} conditions is that the influence of the prior should become negligible as $n$ grows, ensuring that the likelihood dominates. This allows the posterior distribution to be well approximated by an appropriate Gaussian distribution, thereby facilitating the Bernstein–von Mises–type analysis.

With the above notations and notions in hand, we can now state our goal more formally. We aim to show that for large $n$ the posterior distribution of $\uT_1$ (or $\uT_2$) can be well approximated by an appropriate zero mean symmetric matrix variate normal distribution. In other words, we want to show that the total variation norm between $\Pi_{1n}(\cdot\mid\uY^n)$ (or $\Pi_{2n}(\cdot\mid\uY^n)$) and an appropriate zero mean symmetric matrix variate normal distribution will converge in probability to $0$ under $\mathbb{P}_{0}$.

\section{Prior distributions for general covariance or precision matrices}\label{sec3}

Although our main results (Theorem \ref{th_sigma} and \ref{th_omega}) hold for a broad range of prior distributions, it is important to provide specific examples of prior distributions within that class for practical implementation purposes. In this section, we will define some standard and popular prior distributions available for an unstructured covariance or precision matrix. We will show in Section \ref{sec4} that all these prior distributions will satisfy our desired criteria under some mild assumptions.


\subsection{The inverse Wishart prior}
The natural conjugate prior for a covariance matrix is the Inverse Wishart (IW) prior.
We say, $\uSigma\sim IW(\nu+p-1,\;\uPsi_1)$ if the probability density function of $\uSigma$ is given by,
\begin{align}\label{prior1s}
    \pi_{1}^{IW}(\uSigma)\propto \det(\uSigma)^{-(\nu+2p)/2}\exp{(-\tr\left(\uPsi_1\uSigma^{-1}\right)/2)},
\end{align}
where $\nu$ and $\uPsi_1$ are user-specified hyperparameters. It is easy to check the corresponding induced class of the prior distributions on $\uOmega$ will be the Wishart distribution. More precisely if $\uSigma\sim IW(\nu+p-1,\;\uPsi_1)$, then  $\uOmega\sim W(\nu+p-1,\uPsi_1^{-1})$ where
\begin{align}\label{prior1p}
    \pi_{2}^{W}(\uOmega)\propto \det(\uOmega)^{(\nu-2)/2}\exp{(-\tr\left(\uPsi_1\uOmega\right)/2)}.
\end{align}

\noindent While the class of IW priors is a popular choice due to conjugacy and associated algebraic simplicity, it suffers from various drawbacks. \cite{gelman2006prior} strongly discouraged the use of vague inverse gamma priors in a one-dimensional setting, and the IW priors share similar drawbacks in multivariate settings. \cite{alvarez} expounded that a sole degree of freedom parameter regulates the uncertainty for all variance parameters, thereby lacking the flexibility to encompass distinct levels of prior knowledge for various variance components. \cite{Tokuda2011VisualizingDO} discovered that large correlation coefficients correspond to large marginal variances in an IW distribution. This situation can lead to considerable bias in parameter estimations, particularly when correlation coefficients are substantial but marginal variances are limited, and vice versa. Additional comprehensive information can be found in \cite{alvarez, gelman2006prior, gelman2006data, Tokuda2011VisualizingDO}. To overcome these drawbacks several scale mixed versions of IW distributions have been proposed in recent literature. In the subsequent sections, we discuss two prominent members of this class.

\subsection{The diagonal scale-mixed inverse Wishart prior}
The Diagonal Scale-Mixed Inverse Wishart (DSIW) prior for $\uSigma$ is an extension of the Inverse Wishart distribution that incorporates additional parameters to enhance flexibility. Let $\nu > 0$ and $c_{\nu} > 0$ (depending on $\nu$) be user-specified hyperparameters. If $\uDelta$ is a diagonal matrix with the $i$-th diagonal element equal to $\delta_i$, then we can define the DSIW prior using the following hierarchical representation
\begin{align}
\uSigma \mid \uDelta \sim IW(\nu+p-1,\;c_{\nu}\uDelta), \quad \pi(\uDelta) = \prod_{i=1}^p \pi_i(\delta_i),
\end{align}
where $\pi_i(\cdot)$ is a density function with support in the positive real line for every $1 \leq i \leq p$. The marginal prior on $\uSigma$ can be expressed as follows:
\begin{align}\label{prior2s}
\pi_{1}^{DSIW}(\uSigma) \propto \det(\uSigma)^{-(\nu+2p)/2} \prod_{i=1}^p \int_{0}^{\infty} \delta_i^{\nu+p-1/2} \exp{\left\{-\frac{c_{\nu}}{2} (\Sigma^{-1})_{i} \delta_i\right\}} \pi_i(\delta_i) d\delta_i,
\end{align}
where $(\Sigma^{-1})_{i}$ is the $i$-th diagonal element of $\uSigma^{-1}$.
Similarly, the corresponding induced prior on $\uOmega$ is referred to as the diagonal scale mixed Wishart (DSW) prior. It can be expressed as
\begin{align}\label{prior2p}
\pi_{2}^{DSW}(\uOmega) \propto \det(\uOmega)^{(\nu-2)/2} \prod_{i=1}^p \int_{0}^{\infty} \delta_i^{\nu+p-1/2} \exp{\left\{-\frac{c_{\nu}}{2} (\Omega)_{i} \delta_i\right\}} \pi_i(\delta_i) d\delta_i,
\end{align}
\noindent where $(\Omega)_{i}$ is the $i$-th diagonal element of $\uOmega$, and since the Jacobian of the transformation from $\uSigma$ to $\uOmega$ is $\det(\uOmega)^{-(p+1)}$. There are several choices of $\pi$ recommended in the literature. \cite{zava} propose independent log-normal priors on $\delta_i$ with a scale parameter $c_{\nu}=1$ (LN-DSIW prior). Another option they suggest is independent truncated normal priors for the $\delta_i$'s, resulting in the induced prior on $\uSigma$ (TN-DSIW prior). This prior corresponds to the multivariate version of Gelman's folded half-T prior (\cite{gelman2006prior}).

\cite{gelman_book} recommend using independent uniform priors on the $\delta_i$'s with $c_{\nu}=1$ (U-DSIW prior) for non-informative modeling. \cite{huang} suggested employing independent Gamma priors on the $\delta_i$'s with a shape parameter of 2 and $c_{\nu}=2\nu$, (IG-DSIW prior). This prior extends Gelman's Half-t priors on standard deviation parameters to achieve high non-informativity. When $\nu=2$, the correlation parameters under this prior have uniform distributions on the interval $(-1, 1)$. Additionally, for the last two mentioned choices of $\pi_i(.)$'s one can have close form expression for the marginal distribution of $\uSigma$ or $\uOmega$. \cite{gelman2006data} also recommend this prior with $\nu=2$ and $c_{\nu}=1$ to ensure uniform priors on the correlations, similar to the IW prior, but with added flexibility to incorporate prior information about the standard deviations. Similar versions of the aforementioned priors can also be defined for the precision matrix $\uOmega$. See \cite{sarkar} for more detailed information regarding the posterior distributions for these priors in a general framework.

For our analysis, we will later make a mild assumption about the tails of the $\pi_i$s. This assumption holds for all the aforementioned choices of $\pi_i$s, including the LN-DSIW, TN-DSIW, U-DSIW, and IG-DSIW priors. It allows future researchers the flexibility to explore additional options and choose from a wider range of $\{\pi_i\}_{i=1}^p$ distributions.

\subsection{The matrix-$F$ prior}
In the work by \cite{mulder}, a matrix-variate generalization of the $F$ distribution known as the matrix-$F$ distribution for $\uSigma$ is proposed. Similar to the univariate $F$ distribution, the matrix-$F$ distribution can be specified through a hierarchical representation as follows
\begin{align}
\uSigma \mid \bar{\uDelta} \sim IW(\nu+p-1,\;\bar{\uDelta}), \quad \bar{\uDelta} \sim W(\nu^*,\;\uPsi_2),
\end{align}
where $\bar{\uDelta}$ is a matrix-valued random variable, $\nu$ is a positive parameter, $\nu^*$ is the degrees of freedom parameter, and $\uPsi_2$ is a positive definite scale matrix. For the matrix-$F$ distribution, closed-form expressions for the marginal prior on $\uSigma$ and the corresponding induced prior on $\uOmega$ are available. The marginal prior on $\uSigma$ is given by
\begin{align}\label{prior3s}
\pi_{1}^{F}(\uSigma) \propto \det(\uSigma)^{-(\nu+2p)/2} \det(\uSigma^{-1}+\uPsi_2^{-1})^{-(\nu^*+\nu+p-1)/2}.
\end{align}
Similarly, the induced prior on $\uOmega$ can be expressed as
\begin{align}\label{prior3p}
\pi_{2}^{F}(\uOmega) \propto \det(\uOmega)^{(\nu-2)/2} \det(\uOmega+\uPsi_2^{-1})^{-(\nu^*+\nu+p-1)/2}.
\end{align}
The key difference compared to the DSIW prior is that the scale parameter for the base Inverse Wishart distribution is now a general positive definite matrix. For posterior distributions and further details, we refer the reader to \cite{mulder, sarkar}.

\section{BvM results for dense covariance or precision matrices}\label{sec4}
Before providing our main BvM results, we will outline the assumptions made on the true data-generating model and the prior distribution, along with their implications.

\begin{assumption}\label{as1} 
There exists $k_{\sigma}\in(0,1]$ such that
  $\uSigma_{0}\in \mathcal{C_{\sigma}}$, where $\mathcal{C_{\sigma}}=\{\uSigma^{p_n\times p_n}\mid\\\;0<k_{\sigma}\leq\lambda_{min}(\uSigma)\leq\lambda_{max}(\uSigma)$ $\leq1/k_{\sigma}<\infty\}$. We will also assume $\snorm{\uSigma_{0}^{-1/2}Y_1}$ is at most $\sigma_0>0$. Here $k_{\sigma}$ and $\sigma_0$ are fixed constants which do not vary with $n$.
\end{assumption}

\noindent The assumption of uniform boundedness of eigenvalues is a standard assumption in high-dimensional asymptotics for covariance estimation, both in the frequentist and Bayesian settings. It has been widely studied and utilized in various research papers, including \cite{banerjee1, banerjee2, bickel2008regularized, spectrum, xiang}. \cite{bickel2008regularized} referred to the class of covariance matrices satisfying this assumption as ``well-conditioned covariance matrices" and provided several examples of processes that can generate matrices in this class. It is not difficult to check $\uSigma_{0}\in \mathcal{C_{\sigma}}$ iff $\uOmega_{0}\in \mathcal{C_{\sigma}}$.

The bound on the sub-Gaussian norm, involving $\sigma_0$, ensures that there are no unusual or atypical moment behaviors for the distribution of $Y_1$.

\begin{assumption}\label{as2}
 We assume $p_n^5=o(n)$, that is, the number of responses $p_n$ is allowed to grow with $n$, but the ratio $p_n^5/n$ converges to $0$ as $n$ increases. 
\end{assumption}

\noindent As discussed in the introduction, this requirement is unsurprising given the lack of a low-dimensional structure on the covariance matrix and the goal of obtaining BvM results for the \emph{entire} covariance matrix. Relaxing this assumption is challenging without additional structure, such as sparsity, in the covariance or precision matrix. Section \ref{sec7}-\ref{sec9} of our paper is dedicated to demonstrating BvM results under sparsity in the precision matrix. The above assumption can be significantly weakened (see Assumption \ref{asH} in Section \ref{sec9}).

\begin{assumption}\label{as3}
    The sequence of prior distributions $\{\Pi_{1n}(\cdot)\}_{n\geq1}$ (or $\{\Pi_{2n}(\cdot)\}_{n\geq1}$) on $\uSigma$ (or $\uOmega$) is flat around $\uS$ (or $\uS^{-1}$) and the posterior contraction rate under this prior is $\sqrt{\frac{p_n}{n}}$.
\end{assumption}
\noindent When Assumption \ref{as1} holds and $p_n=o(n)$ (which can also be inferred from Assumption \ref{as2}), it can be demonstrated that the sample covariance matrix converges to its true value at a rate of $\sqrt{p_n/n}$ (refer to Lemma 5.2 in \cite{sarkar}). In this setting, it has been shown in \cite{sarkar} that a large class of priors for an unstructured covariance matrix adheres to the contraction rate of $\sqrt{p_n/n}$. Using these results we will show the priors discussed in Section \ref{sec3} (IW, DSIW, and Matrix-$F$) also satisfy this condition under mild assumptions on relevant hyperparameters.

\begin{lemma}\label{postcont1}
   Given Assumptions \ref{as1} and \ref{as2}, the Inverse-Wishart (IW) prior on $\uSigma$, as defined in (\ref{prior1s}), will satisfy Assumption \ref{as3}, even when $\norm{\uPsi_1}=O(p_n)$.
\end{lemma}
\begin{lemma}\label{postcont2}
    Consider a class of DSIW prior distributions on $\uSigma$ as defined in (\ref{prior2s}). Given Assumptions \ref{as1} and \ref{as2}, if there exists a constant $k$ (independent of $n$) such that $\pi_i(x)$ decreases in $x$ for $x > k$ for every $1 \leq i \leq p$, then the prior distribution will satisfy Assumption \ref{as3}.
\end{lemma}
\noindent The proofs of these lemmas are provided in the supplemental document (\cite{supp1}). It is interesting to note that the condition on $\pi_i(.)$ in Lemma \ref{postcont2} is relatively straightforward and encompasses a wide range of commonly used continuous distributions, including truncated normal, half-t distribution, gamma and inverse gamma, beta, Weibull, log-normal, and others. It is worth noting that all the DSIW priors currently discussed in existing literature satisfy this assumption when an appropriate value of $k$ is chosen.

\begin{lemma}\label{postcont3}
   Given Assumptions \ref{as1} and \ref{as2}, the matrix-$F$ prior on $\uSigma$, as defined in (\ref{prior3s}), will satisfy Assumption \ref{as3}, even when $\nu^*=O(p_n)$.
\end{lemma}

\noindent A proof is provided in the supplemental document (\cite{supp1}). With the above lemmas in hand, we now present the key findings of this paper in an unstructured setting.
We first establish the BvM results for the covariance matrix $\uSigma$ in the following theorem. 
\begin{theorem}{\textbf{(BvM Theorem for a Covariance Matrix)}}\label{th_sigma}
 Consider a working Bayesian model for the covariance matrix $\uSigma$ which combines a Gaussian likelihood (\ref{likelihood1}) and utilizes one of the sequences of priors $\{\Pi_{1n}(\cdot)\}_{n\geq1}$ satisfying Assumption \ref{as3}, and assume that true data generating mechanism (see Section \ref{sec2}) satisfies Assumptions \ref{as1}-\ref{as2}. Then
  \begin{align*}
      \int_{B_{1n}} \lvert \pi_{1n}(\uT_1\mid\uY^n)-\phi(\uT_1;\uS)\rvert\; d\uT_1  \overset{P}{\to}0,\; \textit{as}\;n\to\infty\; \textit{under}\;\mathbb{P}_{0},
  \end{align*}
  where $\uT_1=\sqrt{n}(\uSigma-\uS)$ and $\phi(\uT_1;\uS)$ denotes the probability density function of the $\mathcal{SMN}_{p \times p}(\uO,\;2 \uB_p^T\\(\uS\otimes\uS)\uB_p)$ distribution as defined in Section \ref{sec2}. 
\end{theorem}
\noindent Theorem \ref{th_sigma} essentially states that the TV norm between the posterior distribution of $\sqrt{n}(\uSigma-\uS)$ and a $\mathcal{SMN}_{p \times p}(\uO,\;2 \uB_p^T(\uS\otimes\uS)\uB_p)$ converges to zero in probability as $n\to\infty$. In other words, under Assumptions \ref{as1}-\ref{as3}, we can approximate the posterior distribution of $\sqrt{n}(\uSigma-\uS)$ effectively using a symmetric matrix variate normal distribution with mean $\uO$, and with a scale parameter $2 \uB_p^T(\uS\otimes\uS)\uB_p$. This finding proves particularly valuable for constructing credible intervals for (possibly multi-dimensional) functionals of $\uSigma$ directly, as the distribution of $\mathcal{SMN}_{p \times p}(\uO,\;2 \uB_p^T(\uS\otimes\uS)\uB_p)$ can be completely determined from the available data.

We now focus on BvM results for an unstructured precision matrix $\uOmega$ and start with the results for the posterior contraction rate of $\uOmega$.

\begin{lemma}\label{postcont4}
   Suppose the posterior distribution of $\uSigma$ exhibits a contraction rate $\varepsilon_n$, where $\varepsilon_n$ converges to $0$ as $n$ increases. Then under Assumption \ref{as1}, the induced posterior on the precision matrix $\uOmega$ will contract around $\uOmega_0$ at the rate of $\varepsilon_n$ as well.
\end{lemma}

\noindent To prove BvM results for $\uOmega$ with the corresponding induced prior $\Pi_{2n}(\cdot)$, a condition like Assumption \ref{as3} is necessary. Specifically, it is essential for the prior distribution $\Pi_{2n}(\cdot)$ to be flat around $\uS^{-1}$. In fact, under Assumptions \ref{as1} and \ref{as2}, similar results to Lemma \ref{postcont1}, \ref{postcont2}, and \ref{postcont3} can be established for the induced prior on $\uOmega$. The proofs are essentially identical and hence are not provided. We now establish the BvM result for precision matrix $\uOmega$.

\begin{theorem}{\textbf{(BvM Theorem for a Precision Matrix)}}\label{th_omega}
Consider a working Bayesian model which combines a Gaussian likelihood for the precision matrix $\uOmega$ (\ref{likelihood2}) and utilizes one of the sequences of priors $\{\Pi_{2n}(\cdot)\}_{n\geq1}$ satisfying Assumption \ref{as3}, and assume that true data generating mechanism (see Section \ref{sec2}) satisfies Assumptions \ref{as1}-\ref{as2}. Then
  \begin{align*}
      \int_{B_{2n}} \lvert \pi_{2n}(\uT_2\mid\uY^n)-\phi(\uT_2;\uS)\rvert\; d\uT_2  \overset{P}{\to}0,\; \textit{as}\;n\to\infty\; \textit{and under}\;\mathbb{P}_{0},
  \end{align*}
  where $\uT_2=\sqrt{n}(\uOmega-\uS^{-1})$ and  $\phi(\uT_2;\uS)$ denotes probability density function of the $\mathcal{SMN}_{p \times p}(\uO,\;2 \uB_p^T\\(\uS^{-1}\otimes\uS^{-1})\uB_p)$ distribution as defined in Section \ref{sec2}.
\end{theorem}

\noindent The implication of Theorem \ref{th_omega} is the same as Theorem \ref{th_sigma}, except that it applies to $\uOmega$ instead of $\uSigma$. The elements of $\uOmega$ are beneficial when we want to study the conditional dependence structure between the underlying variables.

The proofs for Theorems \ref{th_omega} and \ref{th_sigma} are provided in the Supplementary Material (\cite{supp1}). The key distinction when handling the precision matrix, as opposed to the covariance matrix, lies in the formulation of the likelihood, as illustrated in  (\ref{likelihood1}) and (\ref{likelihood2}). As an expected next step, in the following sections we extend our results from the dense to the sparse setting by introducing the well-known concentration graphical model framework.

\section{Concentration graphical models: preliminaries}\label{sec7}


Before delving into BvM results for concentration graphical models, we will provide the required background material in this section. 
\subsection{Decomposable graphs}
An undirected graph $G = (V, E)$ consists of a vertex set $V = \{1, \ldots, p\}$ with an edge set $E \subseteq \{(i,\;j) \in V \times V : i \neq j\}$, where $(i, j) \in E$ if and only if $(j, i) \in E$. Two vertices $v$ and $v'$ in $V$ are considered adjacent if there exists an edge between them. A complete graph is an undirected graph in which every pair of distinct vertices in $V$ are adjacent. On the other hand, a cycle is a graph that can be represented by a permutation $\{v_1, v_2, \ldots, v_p\}$ of $V$ such that $(v_i, v_j) \in E$ if and only if $\lvert i - j\rvert = 1$ or $\lvert i - j\rvert = p - 1$. The induced subgraph of $G = (V,\;E)$ corresponding to a subset $V' \subseteq V$ is an undirected graph with a vertex set $V'$ and an edge set $E' = E \cap (V' \times V')$. A subset $V'$ of $V$ is considered a clique if the induced subgraph corresponding to $V'$ is a complete graph. Additional details can be found in references such as \cite{letac,lauritzen}. Let $\lvert V \rvert$ denote the cardinality of set $V$. For an undirected graph $G = (V, E)$, we denote $M_G$ as the set of all $\lvert V \rvert \times \lvert V \rvert$ matrices $\uA = (A_{ij})_{1 \leq i,\;j \leq \lvert V \rvert}$ satisfying $A_{ij} = A_{ji} = 0$ for all pairs $(i, j) \notin E$, $i \neq j$. Similarly, $P_G$ represents the set of all symmetric positive definite $(V' \times V')$ matrices that are elements of $M_G$. Now, given the graph $G = (V, E),$ with $V= \{v_1, v_2, \ldots, v_p\}$, we denote $\uA^{>i}=(A_{jk})_{i<j,k\leq p,(i,j)\in E, (i,k)\in E},$ the column vectors $A_{.i}^{>}=(A_{ji})_{j>i, (i,j)\in E}$ and $A_{.i}^{\geq}=(A_{ii}, (A_{.i}^{>})^T)^T$. Also,
\[
A^{\geq i}=
  \begin{bmatrix}
    A_{ii} & (A_{.i}^{>})^T \\
    A_{.i}^{>} & \uA^{>i}
  \end{bmatrix}.
\]
In particular, $A_{.p}^{\geq}=A^{\geq p}=A_{pp}$.
An induced subgraph $G' = (V', E')$ of $G = (V, E)$ is defined when $V' \subseteq V$ and $E' = (V' \times V') \cap E$, and is denoted as $G' \subseteq G$. Let us now revisit the definition of decomposable graphs as stated in \cite{lauritzen}.

\begin{definition}
A graph G is considered decomposable if it does not contain an induced subgraph that forms a cycle of length greater than or equal to $4$.
\end{definition}
\noindent Additional characterizations of decomposable graphs can be found in other references such as \cite{xiang, roverato}. An important property of matrices in the class $P_G$ is worth noting. If $\uOmega\in P_G$, the graph $G$ is decomposable, and the vertices in $V$ are arranged according to a perfect vertex elimination scheme, the Cholesky factor of $\uOmega$ exhibits the same pattern of zeros in its lower triangle (see for example \cite[Theorem $1$]{roverato}).

\subsection{Sparse symmetric matrix-normal distributions}
We introduce a new class of distributions that parallels the symmetric matrix variate normal distributions, called \textit{sparse symmetric matrix-normal distributions (SSMN)}. These distributions will show up as the limiting distributions in the BvM results in Section \ref{sec9}. We start by introducing some useful notations for clarity and convenience. 

Consider a $p \times p$ sparse symmetric matrix $\uA$ where the sparsity structure is given by graph $G$. Let us recall the vectorization of $\uA$ denoted by $\mbox{vec}(\uA)$ as defined in subsection \ref{1.1}. Let $f_p$ represent the number of non-zero unique elements in $\mbox{vec}(\uA)$. Now, consider an $f_p \times 1$ vector $\mbox{vech}^*(\uA)$ that comprises of the non-zero unique elements of $\mbox{vec}(\uA)$. Next, we define a $p^2 \times f_p$ matrix $\uD_{G}$ as an elimination matrix corresponding to graph G (similar to $\uB_p$ mentioned in subsection \ref{1.1}) such that $\mbox{vech}^*(\uA) = \uD_{G}^T\mbox{vec}(\uA)$. We now formally define the SSMN distribution as follows.

\begin{definition}{\textbf{(Sparse Symmetric Matrix-normal Distribution)}}\label{def4}
For a given decomposable graph $G$, let $\uX$ be a $p\times p$ sparse symmetric random matrix taking values in $M_G$. Then $\uX (=\uX^T)$ is said to have a sparse symmetric matrix-normal distribution with parameters 
$\uM \in M_G$, and $\uPsi_1, \uPsi_2$ ($p\times p$ positive definite matrices satisfying $\uPsi_1\uPsi_2=\uPsi_2\uPsi_1$) if $\mbox{vech}^*(\uX)\sim\mathcal{N}_{f_p}(\mbox{vech}^*(\uM),\; \uD_{G}^T(\uPsi_1\otimes\uPsi_2)\uD_{G})$. This distribution is denoted by $\mathcal{SSMN}_{G}(\uM,\;\uD_{G}^T(\uPsi_1\otimes\uPsi_2)\uD_{G})$, and the corresponding probability density function $f(\cdot)$ on $M_G$ is given by 
$$
    f(X)=\frac{\exp\{-\tr(\uPsi_1^{-1}(X-\uM)\uPsi_2^{-1}(X-\uM))/2\}}{(2\pi)^{p(p+1)/4}\det(\uD_{G}^T(\uPsi_1\otimes\uPsi_2)\uD_{G})^{1/2}}. 
$$
\end{definition}

\section{Concentration Graphical Models: model formulation and prior specification for known $G$}\label{sec8}

In a manner similar to Section \ref{sec2}, we consider a set of $n$ independent and identically distributed samples $\uY^{n}=(Y_1,\cdots,Y_n)$ drawn from a multivariate Gaussian distribution $N_p(0,\;\uSigma=\uOmega^{-1})$. For a given undirected graph $G = (V, E)$ with $V = \{1, \ldots, p\}$, the Gaussian concentration model corresponding to G assumes that $\uOmega\in P_G$. Assuming Gaussianity, the log-likelihood function of $\uOmega$, denoted as $l_{3n}(\uOmega)$, can be expressed as follows
\begin{align}\label{likelihood3}
    l_{3n}(\uOmega) = -\frac{np}{2}\log(2\pi) + \frac{n}{2}\log(\det(\uOmega)) - \frac{n}{2}\tr(\uOmega\uS),
\end{align}
where $\uS=\frac{1}{n}\sum_{i=1}^{n}Y_{i}Y_{i}^T$. In a Bayesian framework, we assign a prior $\Pi_{3n}(\cdot)$ to the precision matrix $\uOmega$, with support in $P_G$. For simplicity of notation, we will sometimes refer to the prior density $\pi_{3n}(\cdot)$ as $\pi_{3}(\cdot)$.

We will now specify the \textit{true data-generating mechanism} within the above framework. We assume that the observations $Y_1,\dots, Y_n$ are independently and identically distributed from a multivariate Gaussian distribution $N_{p_n}(0,\;\bar{\uOmega}_n^{-1})$, where $\{{\bar{\uOmega}_n}\}_{n\geq1}$ represents the sequence of true precision matrices. Let $G_n=(V_n,\; E_n)$, with $V_n=\{1,\cdots,p_n\}$, be a decomposable graph where the vertices are ordered according to a perfect vertex elimination scheme. We will assume that $\bar{\uOmega}_n\in P_{G_n}$. Let $d_n$ denote the maximum number of non-zero entries in any row of the symmetric matrix $\bar{\uOmega}_n$. Also, define $a^{G}(=a^{G_n})$ is the product of $\max_{1\leq j \leq p_n} n_j +1$ and  $\max_{1\leq i \leq p_n} r_i +1$. Here $n_j=\{i:1\leq j <i \leq p, (i,j) \in E_n\}$ and $r_i=\{j:1\leq j < i \leq p, (i,j) \in E_n\}$. We denote the probability measure underlying the true model as $\mathbb{P}_{0,\; G_{n}}$. For simplicity, we will use $\mathbb{P}_{0,G}$ instead of $\mathbb{P}_{0,\;G_{n}}$. 
Next, we define the maximum likelihood estimator of $\uOmega$ within the class $P_{G_n}$ as
\begin{align}\label{omega_g}
    \hat{\uOmega}_G(\;=  \hat{\uOmega}_{G_n})=\underset{\Omega\in P_{G_n}}{\sup} l_{3n}(\uOmega),
\end{align}
where $l_{3n}(\uOmega)$ is defined in  (\ref{likelihood3}).Let $\uT_3=\sqrt{n}(\uOmega-\hat{\uOmega}_G)$ be a centered version of $\uOmega$. In this context, we define the function
\begin{align}
M_{3n}(\uT_3)=\exp{\left(l_{3n}\left(\hat{\uOmega}_G+\dfrac{\uT_3}{\sqrt{n}}\right)-l_{3n}\left(\hat{\uOmega}_G\right)\right)},\label{m3}
\end{align}
where $\uT_3$ belongs to $B_{3n}$, and $B_{3n}=\{\uT_3:\hat{\uOmega}_G+\frac{\uT_3}{\sqrt{n}}\in P_{G_n}\}$. If $\uT_3$ falls outside $B_{3n}$, we set $M_{3n}(\uT_3)$ to be zero. Clearly $B_{3n}$ is a subset of $M_{G_n}$. Now, suppose the posterior distribution for $\uT_3$ is given by $\Pi_{3n}(\cdot\mid\uY^n)$. Analogously, let $\pi_{3n}(\cdot\mid\uY^n)$ represent the corresponding posterior density. Then it is not difficult to check that,
\begin{align}
    &\pi_{3n}(\uT_3\mid\uY^n)=\frac{M_{3n}(\uT_3)\pi_{3}\left(\hat{\uOmega}_G+\frac{\uT_3}{\sqrt{n}}\right)}{\bigintsss_{B_{3n}} M_{3n}(W)\pi_{3}\left(\hat{\uOmega}_G+\frac{W}{\sqrt{n}}\right))dW}\label{post3}.
\end{align}
Our objective is to demonstrate that the total variation norm between $\Pi_{3n}(\cdot\mid\uY^n)$ and an appropriate zero-mean sparse symmetric matrix variate normal distribution converges in probability to $0$ under $\mathbb{P}_{0,G}$.

\subsection{The $G$-Wishart distribution}
As in Section \ref{sec3} while our main results hold for a broad range of prior distributions, it is important to provide specific examples of prior distributions within that class for practical implementation purposes. In this subsection, we will define a standard prior distribution available for the precision matrix under the Gaussian concentration model concerning the graph $G$ defined in Section \ref{sec7}.
\cite{dawid} developed a class of hyper inverse Wishart distributions for $\uSigma=\uOmega^{-1}$ when $\uOmega\in P_G$. The corresponding class of induced priors for $\uOmega$ is known as the class of $G$-Wishart distributions on $P_G$ (See \cite{atay, roverato}). Specifically, the $G$-Wishart distribution with parameters $\beta\geq 0$ and $\uPsi_3$ positive definite, denoted by $W_G(\beta,\uPsi_3)$, has a density proportional to
\begin{align}\label{prior4p}
    \pi_{3}^{WG}(\uOmega)\propto \det(\uOmega)^{\beta/2}\exp{(-\tr\left(\uPsi_3\uOmega\right)/2)},\;\uOmega\in P_G.
\end{align}
The class of $G$-Wishart distributions on $P_G$ forms a conjugate family of priors under the Gaussian concentration graphical model corresponding to $G$. If $G$ is decomposable, then quantities such as the mean, mode, and normalizing constant for $W_G(\beta,\uPsi_3)$ are available in closed form (see, for instance, \cite{rajaratnam}). In Section \ref{sec9}, we will demonstrate that the $G$-Wishart distribution falls into our desired class of prior distributions under suitable assumptions.

\section{BvM and posterior consistency results for sparse precision matrices in the known-$G$ Case}\label{sec9}

As mentioned earlier in Section \ref{sec4}, achieving BvM results hinges on the rates of posterior contraction. As mentioned in the introduction, rates in the existing literature either lack optimality or relies on stronger assumptions to establish contraction rates for the precision matrix within the Gaussian concentration model framework. Therefore, before presenting our main BvM results, we will refine posterior contraction rates for the precision matrix within this framework at least when the underlying graph is known. Initially, we will outline a set of standard assumptions (which are more relaxed compared to previous work in \cite{banerjee1, banerjee2, Lee, liu2019empirical, xiang}) that are necessary to achieve these posterior contraction rates, along with their brief implications.
\renewcommand{\theassumption}{\Alph{assumption}}
\makeatletter
\makeatother

\begin{assumption}\label{asD} 
The eigenvalues of $\{{\bar{\uOmega}_n}\}_{n\geq1}$ are uniformly bounded i.e. $\bar{\uOmega}_n\in \mathcal{C_{\sigma}}$, where $\mathcal{C_{\sigma}}$ is defined in Assumption \ref{as1}.
\end{assumption}
\noindent As noted in Section \ref{sec4}, the assumption of uniformly bounded eigenvalues is standard in high-dimensional asymptotics for precision matrix estimation. This assumption aligns with those employed in the Bayesian framework within the Gaussian concentration model literature, as evidenced by \cite{banerjee1, banerjee2, Lee, liu2019empirical, Yabo_Niu, xiang}. 
The next condition is needed for establishing a posterior contraction rate in the spectral norm.
\begin{assumption}\label{asE}
    $a^{G_n}\log p_n=o(n)$ and $p_n\to\infty$ as $n\to\infty$.
\end{assumption}
If the Frobenius norm is utilized instead of the spectral norm, Assumption \ref{asE} requires adjustment as the Frobenius norm is larger in magnitude thus resulting in the following assumption.
\begin{assumption}\label{asF}
    $(p_n+ |E_n|) \log p_n=o(n)$ and $p_n\to\infty$ as $n\to\infty$.
\end{assumption}
The final assumption imposes mild restrictions on the hyper-parameters corresponding to the $G$-Wishart prior distribution.

\begin{assumption}\label{asG}
For each $n \ge 1$, we place the prior $W_{G_n}(\beta, \uPsi_{3n})$ on the concentration matrix $\uOmega$, where $\beta > 0$ is fixed. The eigenvalues of $\uPsi_{3n}$ are uniformly bounded, that is, $\uPsi_{3n} \in \mathcal{C}_{\sigma^*}$ for some constant $\sigma^*$, where $\mathcal{C}_{\sigma^*}$ is defined in Assumption \ref{as1}.
\end{assumption}

\noindent To achieve a posterior contraction rate in the spectral norm, \cite{banerjee1,xiang} imposes the assumption $d_n^4\log p_n=o(n)$ for spectral norm consistency, and $d_n^5\log p_n =o(n)$ for matrix $(\infty,\infty)$-norm consistency. In \cite{Lee}, the authors refine the matrix $(\infty,\infty)$-norm consistency constraint to $d_n^4\log p_n=o(n)$ in a setting where the underlying graph/sparsity structure $G$ is assumed to be unknown. Assumption \ref{asE} represents a significantly more lenient assumption within the Bayesian framework compared to those utilized in the existing literature. For a simple demonstration, consider a star graph formed with $p_n$ variables or nodes. This graph is essentially a tree with $p_n$ nodes, where one node has a vertex degree of $p_n-1$ and the remaining $p_n-1$ nodes have a vertex degree of $1$. It is evident that graph is decomposable since it does not contain any cycle. In this simple scenario, Assumption \ref{asE} translates to the optimal condition $p_n\log p_n=o(n)$, whereas \cite{banerjee1,xiang} necessitates $p_n^4\log p_n=o(n)$. In simpler terms, Assumption \ref{asE} is considerably more relaxed compared to the assumptions made in \cite{banerjee1,xiang,Lee} when both $d_n$ and $a^{G_n}$ are of the order $p_n$.

When considering Frobenius norm convergence rates, \cite{liu2019empirical} require Assumption \ref{asF} in a setting where the underlying graph/sparsity structure $G$ is not known, but an additional assumption that $p_n= O(n^{\delta})$ for some $\delta \in (0,1)$ is imposed. However, we are able to get rid of this additional assumption in the currently known $G$ setting, rendering our approach more applicable to high-dimensional situations.

We now state our main posterior consistency results (Theorem \ref{postthm1} \& Theorem \ref{postthm2}) for the graphical model setting. The proofs are available in the supplemental document (\cite{supp1}).
\begin{theorem}{\textbf{(Posterior Contraction Rate for a Sparse Precision Matrix under spectral norm)}}\label{postthm1}
Let $Y_1, Y_2, Y_3, \cdots, Y_n$ be independent and identically distributed Gaussian random variables with mean $0$ and precision matrix $\uOmega\in P_{G_n}$, where $G_n$ is the undirected graph encoding the sparsity in $\uOmega$. Consider a working Bayesian model that utilizes a sequence of $G$-Wishart priors satisfying Assumption \ref{asG}, and assume that true data generating mechanism (see Section \ref{sec8}) satisfies Assumptions \ref{asD} and \ref{asE}. Then
  \begin{align*}
    \Pi_{3n}\left(\norm{\uOmega-\bar{\uOmega}} > M \sqrt{\frac{a^{G_n}\log p_n}{n}} \mid \uYn\right) \overset{P}{\to} 0
  \end{align*}
  as $n \to \infty$, under $\mathbb{P}_{0,G}$ and for sufficiently large $M$.
\end{theorem}

\noindent Similar to our previous discussion, this contraction rate is notably more lenient in many scenarios compared to the contraction rate of $\sqrt{d_n^4\log p_n /n}$ (as seen in \cite{banerjee1,xiang}). For instance, in the case of the star graph, our contraction rate simplifies to $\sqrt{p_n\log p_n /n}$, whereas it becomes $\sqrt{p_n^4\log p_n /n}$ for \cite{banerjee1,xiang}. If one opts for using the matrix $(\infty,\infty)$ norm instead of the spectral norm in this scenario, our contraction rate is equal to $\sqrt{p_n^2\log(p_n)/n}$, significantly surpassing the $\sqrt{p_n^4\log(p_n)/n}$ obtained by \cite{Lee}. When $d_n$ remains constant or increases slowly with $n$ (as observed in scenarios such as the banded concentration graphical model from \cite{banerjee1}), all these contraction rates are relatively comparable.

\cite[Theorem~5]{cai_minimax} showed that if the true precision matrix $\bar{\uOmega}$ is restricted to a class in which each nonzero entry is bounded below by $\sqrt{\log p_n / n}$, and $\bar{\uOmega}$ has uniformly bounded eigenvalues (Assumption~\ref{as1}), the minimax risk for estimating a sparse precision matrix under the spectral norm is of order $\sqrt{d_n^{2}\log p_n / n}$, provided that $d_n^{2} (\log p_n)^{3} / n = O(1)$. 
Since $a_{G_n}$ is bounded above by $d_n^{2}$, this implies that in the worst-case scenario, under the stronger condition $d_n^{2} (\log p_n)^{3} / n = O(1)$, our posterior contraction rate coincides with the frequentist minimax rate. It is worth noting, however, that our parameter space is restricted to decomposable graphs, whereas the bounds in \cite{cai_minimax} apply to a more general class of sparse precision matrices.


In the next theorem, we will describe the posterior contraction rate under the Frobenius norm.
\begin{theorem}{\textbf{(Posterior Contraction Rate for a Sparse Precision Matrix Under the Frobenius Norm)}}\label{postthm2}
Let $Y_1, Y_2, Y_3, \cdots, Y_n$ be independent and identically distributed Gaussian random variables with mean $0$ and precision matrix $\uOmega\in P_{G_n}$, where $G_n$ is the undirected graph encoding the sparsity in $\uOmega$. Consider a working Bayesian model that utilizes a sequence of $G$-Wishart priors satisfying Assumption \ref{asG}, and assume that true data generating mechanism (see Section \ref{sec8}) satisfies Assumptions \ref{asD} and \ref{asF}. Then
  \begin{align*}
    \Pi_{3n}\left(\fnorm{\uOmega-\bar{\uOmega}} > M \sqrt{\frac{(p_n+ |E_n|)\log p_n}{n}} \mid \uYn\right) \overset{P}{\to} 0
  \end{align*}
  as $n \to \infty$, under $\mathbb{P}_{0,G}$ and for sufficiently large $M$.
\end{theorem}
\noindent This contraction rate aligns with the optimal contraction rate for maximum likelihood estimators of a sparse precision matrix in the frequentist setup (see \cite{rothman, Lam}). In the Bayesian paradigm, \cite{liu2019empirical} accomplish a similar contraction rate, As discussed previously, \cite{liu2019empirical} consider the case where $G$ is unknown, and need an additional assumption $p_n= O(n^{\delta})$ for some $\delta \in (0,1)$ to achieve such a contraction rate. Our results show that the additional assumption is unnecessary in the known $G$ setting.

With the posterior contraction results in hand, we now present our BvM theorem for the sparse precision matrix within the Gaussian graphical model framework. 
It is important to acknowledge that proving BvM results often demands stronger assumptions than those required for posterior consistency. An analogy can be drawn to the frequentist setup, where demonstrating results akin to the Central Limit Theorem typically calls for stronger assumptions compared to those needed for establishing simple parameter consistency. With this in mind, we state the two regularity assumptions required for our BvM result. 
\begin{assumption}\label{asH}
    $\min(p_n^2(a^{G_n})^3, (p+|E_n|)^3)=o(n/(\log p_n)^3)$ and $p_n\to\infty$ as $n\to\infty$.
\end{assumption}
\noindent Note that $a^{G_n}$ is bounded above by $d_n^{2}$. Here $d_n$, which represents the maximum number of non-zero entries in any row of the true precision matrix $\bar{\uOmega}$, stays constant or increases slowly with $n$ (e.g. as observed in cases such as the banded concentration graphical model from \cite{banerjee1}), Assumption \ref{asH} becomes much weaker compared to Assumption \ref{as2} ($p_n^5=o(n)$). Also, for the example involving the star graph (see the discussion after Assumption \ref{asG}), Assumption \ref{asH} simplifies to $(p_n\log p_n)^3=o(n)$.
This implies, as expected, that we can relax the strict conditions on $p_n$ (for BvM in the unstructured setting) by imposing a sparse structure on the precision matrix and handling scenarios with larger $p_n$ without sacrificing the validity of our results.
Similar to Assumption \ref{as3} in the unstructured case, we need a flatness assumption on the prior distribution for the sparse precision matrix. Let $\varepsilon_{3,n}= \sqrt{a^{G_n}\log p_n/n}$ be the posterior contraction rate in Theorem \ref{postthm1} and define $\rho^{\pi_{3}}(\varepsilon_{3,n})$  as
\begin{align}\label{flat}
\rho^{\pi_{3}}(\varepsilon_{3,n}):= \underset{\mathbf{T_3}\in D(\varepsilon_{3,n})}{\sup}\;\left\lvert\frac{\pi_{3}(\hat{\uOmega}_{G_n}+\dfrac{\uT_3}{\sqrt{n}})}{\pi_{3}(\hat{\uOmega}_{G_n})}-1\right\rvert,
\end{align}
where $\pi_{3}(\cdot)$ denotes the prior density for $\uOmega$, $B_{3n}=\{\uT_3:\hat{\uOmega}_G+\frac{\uT_3}{\sqrt{n}}\in P_{G}\}$ and $D(\varepsilon_{3,n})=\{\uT_3\in B_{3n}\mid\norm{\uT_3} \leq \sqrt{n}\varepsilon_{3,n} \}$. Now, we will formally state the flatness assumption for the sequence of prior distributions ${\Pi_{3n}(\cdot)}_{n\geq1}$.
\begin{assumption}\label{asI}
    The sequence of prior distributions $\{\Pi_{3n}(\cdot)\}_{n\geq1}$ on $\uOmega$ with support on $P_{G_n}$ is flat around $\hat{\uOmega}_{G_n}$ i.e. $\rho^{\pi_{3}}(\varepsilon_{3,n})\to 0$ in probability as $n\to\infty$.
\end{assumption}
\noindent  When Assumption \ref{asD} and \ref{asE} are satisfied it has been shown that the maximum likelihood estimator, $\hat{\uOmega}_G$ converges at the rate $\varepsilon_{3,n}$ (\cite{xiang,KRRZ:2024}). The lemma below uses this to show that the $G$-Wishart prior distribution satisfies Assumption \ref{asI} under certain conditions on the hyperparameters. A proof is provided in the supplemental document (\cite{supp1}).
\begin{lemma}\label{gw}
Let $W_{G_n}(\beta,\uPsi_{3,n})$ denote the $G$-Wishart prior distribution defined in (\ref{prior4p}). Under Assumptions \ref{asD}, \ref{asH}, and \ref{asG}, this prior distribution satisfies Assumption \ref{asI}.
\end{lemma}

\noindent We now present the key result of this section.
\begin{theorem}{\textbf{(BvM Theorem for a Sparse Precision Matrix with Known Sparsity)}}\label{th_omega_sp}
Let $Y_1, Y_2,\\Y_3, \cdots, Y_n$ be independent and identically distributed Gaussian random vectors with mean $0$ and precision matrix $\uOmega$, where $\uOmega\in P_{G_n}$. Here $G_n$ is the undirected graph which encodes the sparsity in $\uOmega$. Consider a working Bayesian model that utilizes one of the sequences of priors $\{\Pi_{3n}(\cdot)\}_{n\geq1}$ satisfying Assumption \ref{asI}, and assume that true data generating mechanism (see Section \ref{sec8}) satisfies Assumptions \ref{asD} and \ref{asH}, Then
  \begin{align*}
      \int_{B_{3n}} \lvert \pi_{3n}(\uT_3\mid\uY^n)-\phi(\uT_3;\hat{\uOmega}_G)\rvert\; d\uT_3 \overset{P}{\to}0,\; \textit{as}\;n\to\infty 
  \end{align*}
  where $\uT_3=\sqrt{n}(\uOmega-\hat{\uOmega}_G)$ and $\phi(\uT_3;\hat{\uOmega}_G)$ denotes the probability density function of the $\mathcal{SSMN}_{G_n}(\uO,\\\;2 \uD_{G}^T(\hat{\uOmega}_G\otimes\hat{\uOmega}_G)\uD_{G})$ distribution as defined in Section \ref{sec7}.
\end{theorem}
\noindent The proof of Theorem \ref{th_omega_sp} is included in the supplemental document (\cite{supp1}).

\section{Extension of BvM results when the underlying graph is unknown}\label{sec10}

In Section \ref{sec9}, we assumed that the true decomposable graph $G$, which encodes the sparsity structure in the precision matrix $\uOmega$, is completely known. However, this is rarely the practical case. Typically, the underlying sparsity structure in $\uOmega$ is unknown. In this section, we will demonstrate that even when $G$ is unknown, we can still establish a BvM result similar to Theorem \ref{th_omega_sp}, provided we have a mechanism to consistently estimate $G$. In the process, we will also establish posterior contraction rates for $\uOmega$ in the unknown $G$ setting. These rates are shown to match the contraction rates for the known graph setting (Theorems \ref{postthm1} and \ref{postthm2}) under suitable conditions.
\subsection{Hirearchical $G-$Wishart prior}

Since the graph $G=(V, E)$ is unknown, we also need to specify a prior distribution on the graph $G$. In particular, we consider a hierarchical $G-$Wishart prior on $(G,\uOmega)$ given by 
\begin{align}\label{hir_G}
    &\uOmega\mid G\sim W_G(\beta,\uPsi_3),\;\uOmega \in P_{G}\notag\\
    & G\sim \pi(G),\;G \in\mathcal{D},
\end{align}
where $\beta$ and $\uPsi_3$ are the corresponding hyperparameters of the $G-$Wishart distribution as described in (\ref{prior4p}) and $\mathcal{D}$ is a set of all decomposable graphs. A very popular choice of $\pi(G)$ found in the literature is given by
\begin{align}\label{G_prior}
    \pi(G)\propto {p(p-1)/2 \choose |E|}^{-1}\exp\{-|E|\tau\log p\}\;\mathbbm{1}\{G\in \mathcal{D},|E|\leq R\},
\end{align}
for some constant $\tau>0$ and a positive integer $R$. The condition $|E| \leq R$ implies that we focus only on graphs that do not have too many edges. Under the prior $\pi(G)$, the prior mass decreases exponentially with respect to the graph size $|E|$, and given a graph size, the locations of edges are sampled from a uniform distribution. These priors have been employed both in the high-dimensional regression and covariance estimation literature, see for example \cite{castillo1,yang,Martin,liu2019empirical, Lee,banerjee2,cao}. Notably, if we choose $\tau = \log\left(\frac{1-q}{q}\right)$, then for a fixed graph size $|E|$, this prior essentially prior collapses to the classical Erdős–Rényi (Bernoulli‐graph) prior form
\begin{align}
    \pi(G) \propto q^{|E|}(1-q)^{{p \choose 2} - |E|} \;\mathbbm{1}\{G \in \mathcal{D}, |E| \leq R\},
\end{align}
which means we are placing an i.i.d. Bernoulli prior (usual over all the edges of the graph with success probability $q$.

However, in our analysis, we will not rely on any specific form of $\pi(G)$. Instead, we will impose a generic condition (see Assumption \ref{asJ}) that is satisfied by the prior distribution described in (\ref{G_prior}). This approach allows us to unify the results in a more general setting and leaves room for applicability to other priors in the future.

\subsection{Model formulation and posterior distributions for Unknown $G$}\label{subsec9.2}

Similar to Section \ref{sec8}, we consider $n$ independent and identically distributed samples $\mathbf{Y}^{n} = (Y_1, \ldots, Y_n)$ drawn from a multivariate Gaussian distribution $N_p(0, \boldsymbol{\Sigma} = \boldsymbol{\Omega}^{-1})$. For a given undirected graph $G = (V, E)$ with $V = \{1, \ldots, p\}$, the Gaussian concentration model corresponding to $G$ assumes that $\boldsymbol{\Omega} \in \mathcal{P}_G$. In a Bayesian framework, we assign a prior $\Pi_{4n}(\cdot \mid G)$ on $\boldsymbol{\Omega}$ given the graph $G$ with support in $\mathcal{P}_G$, and a prior $\Pi_n(G)$ on the graph $G$. For notational simplicity, we will sometimes refer to the prior density $\pi_{4n}(\cdot\mid G)$ as $\pi_{4}(\cdot\mid G)$ and $\pi_n(G)$ as $\pi(G)$.

We now specify the \textit{true data-generating mechanism.} We assume that the observations $Y_1,\dots, Y_n$ are independently and identically distributed from a multivariate Gaussian distribution $N_{p_n}(0,\;\bar{\uOmega}_n^{-1})$, where $\{{\bar{\uOmega}_n}\}_{n\geq1}$ represents the sequence of true precision matrices. Let $G_{0n}=(V_{0n},\; E_{0n})$, with $V_{0n}=\{1,\cdots,p_n\}$, be a decomposable graph encoding the sparsity in $\bar{\uOmega}_n$, i.e., $\bar{\uOmega}_n\in P_{G_{0n}}$. We assume that the vertices of $G_{0n}$ are ordered according to a perfect vertex elimination scheme. Let $d_n$ denote the maximum number of non-zero entries in any row of the symmetric matrix $\bar{\uOmega}_n$. Also, define $a^{G_{0}}(=a^{G_{0n}})$ as the product of $\max_{1\leq j \leq p_n} n_j +1$ and  $\max_{1\leq i \leq p_n} r_i +1$. Here $n_j=\{i:1\leq j <i \leq p, (i,j) \in E_n\}$ and $r_i=\{j:1\leq j < i \leq p, (i,j) \in E_{0n}\}$. We denote the probability measure underlying the true model as $\mathbb{P}_{0,\; G_{0n}}$. For simplicity, we will use $\mathbb{P}_{0,G_{0}}$ instead of $\mathbb{P}_{0,\;G_{0n}}$.

When the graph $G$ is known, as discussed in Section \ref{sec8}, we center our precision matrix using the corresponding maximum likelihood estimate $\hat{\uOmega}_G$, enabling us to establish BvM results. However, in the case where $G$ is unknown, this estimate is unavailable. We can overcome this issue by employing a two-stage estimation technique to construct another likelihood-based estimate of $\uOmega$. Note that several methods for finding consistent estimates of $G$ in high-dimensional settings are available in the literature, see for example  \cite[Theorem 2]{raskutti2008model}, \cite[Theorem 2]{khare}. First, we use one such method to obtain a consistent estimator $\hat{G}_n$ of $G$ (under $\mathbb{P}_{0,\;G_{0}}$), and then, given $\hat{G}_n$, we calculate the conditional maximum likelihood estimate of $\boldsymbol{\Omega}$ as in (\ref{omega_g}). In particular, we use the estimator 
\begin{align}\label{omega_g1}
    \hat{\uOmega}_{\hat{G}}(\;=  \hat{\uOmega}_{\hat{G}_n })=\underset{\Omega\in P_{\hat{G}_n }}{\sup} l_{3n}(\uOmega),
\end{align}
where $l_{3n}(\uOmega)$ is defined in (\ref{likelihood3}), for centering purposes and follow a similar process as in Section \ref{sec8}. Let $\uT_4=\sqrt{n}(\uOmega-\hat{\uOmega}_{\hat{G}})$ be a centered version of $\uOmega$. In this context, we define the function
\begin{align}
M_{4n}(\uT_4\mid G)=\exp{\left(l_{3n}\left( \hat{\uOmega}_{\hat{G}}+\dfrac{\uT_4}{\sqrt{n}}\right)-l_{3n}\left(\hat{\uOmega}_{\hat{G}}\right)\right)},\label{m3}
\end{align}
where $\uT_4$ belongs to $B_{4n}$, and $B_{4n}=\{\uT_4:\hat{\uOmega}_{\hat{G}}+\frac{\uT_4}{\sqrt{n}}\in P_{G}\}$. If $\uT_4$ falls outside $B_{4n}$, we set $M_{4n}(\uT_4\mid G)$ to be zero. A notable difference compared to the known $G$ setting is that $B_{4n}$ is a subset of the space of all possible $p\times p$ real matrices and not of $M_G$ as in Section \ref{sec8}. Let the marginal posterior distribution for $\uT_4$ and $G$ be denoted by $\Pi_{4n}(\uT_4\mid\uY^n)$ and $\Pi_{n}(\cdot\mid\uY^n)$ respectively. Analogously, let $\pi_{4n}(\cdot\mid\uY^n)$ and $\pi_{n}(\cdot\mid\uY^n)$ represent the corresponding posterior densities. Then it follows that 
\begin{align}
   \pi_{4n}(\uT_4\mid\uY^n) = \sum_{ G \in \mathcal{D}} \pi_{4n}(\uT_4\mid\uY^n, G) \pi_{n}(G\mid\uY^n),
\end{align}
where
\begin{align}
    &\pi_{4n}(\uT_4\mid\uY^n, G)=\frac{M_{4n}(\uT_4)\pi_{4}\left(\hat{\uOmega}_{\hat{G}}+\frac{\uT_4}{\sqrt{n}}\mid G\right)}{\bigintsss_{B_{4n}} M_{4n}(W)\pi_{4}\left(\hat{\uOmega}_{\hat{G}}+\frac{W}{\sqrt{n}}\mid G\right)dW}\label{post3}.
\end{align}
We aim to show that, under $ \mathbb{P}_{0, G_0} $, the total variation norm between $ \Pi_{4n}(\cdot \mid \mathbf{Y}^n) $ and a suitable zero-mean sparse symmetric matrix variate normal distribution converges in probability to $0$. The primary challenge lies in the fact that now $ \mathbf{T}_4$ is a real matrix of size $ p \times p $, and we lack information regarding its sparsity structure. This significantly complicates the analysis compared to the known graph scenario in Section \ref{sec8}.

\section{BvM and posterior consistency results for sparse precision matrices in the unknown $G$ case}\label{sec9}

In this section, we first establish analogs of Theorems \ref{postthm1} and \ref{postthm2} for an unknown graph $G$, and then we proceed to derive a BvM result analogous to Theorem \ref{th_omega_sp}. However, to establish a posterior contraction rate when the underlying sparsity structure encoded by $ G $ is unknown, we require a property known as \textit{strong graph selection consistency}. In the covariance estimation literature, strong graph selection consistency is often used as a stepping stone to establish the Bernstein-von Mises theorem in settings where the true graph $G$ is unknown. See, for example, \cite{FangGhosh2024, RyanNing2020}. We impose this property through the following assumption. 
\begin{assumption}\label{asJ}
    $ \pi_{n}(G_0 \mid \mathbf{Y}^n) \overset{P}{\to} 1 $ as $ n \to \infty $ under $ \mathbb{P}_{0,G_{0}} $.
\end{assumption}
\noindent \cite{Lee} established strong selection consistency for the prior structure described in (\ref{hir_G}) and (\ref{G_prior}), demonstrating this result under a set of standard regularity assumptions on the true data-generating model. We briefly summarize the high-level ideas behind these assumptions below.
First, they assumed the growth condition $|E_{0n}| \leq R$, where $R$ is the prior cut-off value (see (\ref{G_prior})), ensuring that the true graph $G_0$ is not too large and therefore lies within the prior support. A similar assumption is also required in our posterior contraction result (see Theorem~\ref{postthm5}). The next two assumptions concern upper and lower bounds on relevant partial correlations of the true covariance matrix $\bar{\uSigma}_n = \bar{\uOmega}_n^{-1}$, expressed in terms of suitable functions of $n$ and $p_n$. One of these assumptions is intended to rule out the possibility that there exists a set of variables $S$ with $|S| \leq 3R$ such that two variables (connected by an edge in $G_0$) become perfectly linearly dependent when conditioning on $S$. The other assumption is the minimum signal size condition (the well-known beta-min assumption). As discussed prior to Theorem~\ref{postthm5}, this type of condition is essential for establishing strong selection consistency but is not required when the goal is only posterior consistency. Regarding the behavior of $R$, \cite{Lee} assumed $R \propto \left( \frac{n}{\log(\max(p_n,n))} \right)^{\xi/3}
\quad \text{for some } \xi \in [0,1]$, which is more restrictive than what is required when one is interested solely in posterior consistency (see Theorem~\ref{postthm5}). For further technical details, we refer the reader to \cite[Section~3]{Lee}.
We now state our posterior contraction results (Theorem \ref{postthm3} \& Theorem \ref{postthm4}) for the Gaussian graphical models setting with unknown graph structure. The proofs are available in the supplemental document (\cite{supp1}).
\begin{theorem}{\textbf{(Posterior Contraction Rate for a Sparse Precision Matrix Under the Spectral Norm)}}\label{postthm3}
Let $Y_1, Y_2, Y_3,\cdots, Y_n$ be independent and identically distributed Gaussian random variables with mean $0$ and precision matrix $\uOmega\in P_{G}$. Consider a working Bayesian model that utilizes a hierarchical $G$-Wishart prior on $(\uOmega, G)$ as described in (\ref{hir_G}), satisfying Assumption \ref{asG} and \ref{asJ}, and assume that true data generating mechanism (see subsection \ref{subsec9.2}) satisfies Assumptions \ref{asD} and \ref{asE}. Then
  \begin{align*}
    \Pi_{4n}\left(\norm{\uOmega-\bar{\uOmega}} > M \sqrt{\frac{a^{G_{0n}}\log p_n}{n}} \mid \uYn\right) \overset{P}{\to} 0
  \end{align*}
  as $n \to \infty$, under $\mathbb{P}_{0,G_0}$ and for a sufficiently large constant $M$.
\end{theorem}

In the next theorem, we establish posterior contraction rates under the Frobenius norm.
\begin{theorem}{\textbf{(Posterior Contraction Rate for a Sparse Precision Matrix Under the Frobenius Norm)}}\label{postthm4}
Let $Y_1, Y_2, Y_3, \cdots, Y_n$ be independent and identically distributed Gaussian random variables with mean $0$ and precision matrix $\uOmega\in P_{G}$. Consider a working Bayesian model that utilizes a sequence of hierarchical $G$-Wishart prior on $(\uOmega, G)$ as described in (\ref{hir_G}), satisfying Assumption \ref{asG} and \ref{asJ}, and assume that true data generating mechanism (see subsection \ref{subsec9.2}) satisfies Assumptions \ref{asD} and \ref{asF}. Then
  \begin{align*}
    \Pi_{4n}\left(\fnorm{\uOmega-\bar{\uOmega}} > M \sqrt{\frac{(p_n+ |E_{0n}|)\log p_n}{n}} \mid \uYn\right) \overset{P}{\to} 0
  \end{align*}
  as $n \to \infty$, under $\mathbb{P}_{0,G}$ and for sufficiently large constant $M$.
\end{theorem}

Thus, we can obtain precisely the same contraction rate for the precision matrix $\uOmega$ even when the underlying graph structure is unknown, assuming we have a strong selection consistency result. The implications of Theorems \ref{postthm3} and \ref{postthm4} will mirror those of Theorems \ref{postthm1} and \ref{postthm2}, as discussed in Section \ref{sec9}.

\begin{remark}
Prompted by a reviewer’s suggestion, we examined whether posterior contraction rates can be established without Assumption~\ref{asJ}, that is, without assuming strong selection consistency. Strong selection consistency typically relies on a minimum signal size, or the well-known \emph{beta-min condition}. While such conditions are often essential for exact sparsity recovery, they are not generally required when the goal is only to establish posterior contraction rates. In contrast, when proving Bernstein-von Mises-type results (which, in our setting, build on these contraction rates), strong selection consistency becomes difficult to avoid (see, e.g., \cite{RyanNing2020,FangGhosh2024}). However, if one is primarily interested in estimation accuracy, such a condition is not often required. The next theorem states the posterior contraction rate in the Frobenius norm that holds without assuming strong selection consistency.
\end{remark}
\begin{theorem}{\textbf{(Posterior Contraction Rate for a Sparse Precision Matrix without minimum signal size assumption)}}\label{postthm5}
Let $Y_1, Y_2, Y_3, \cdots, Y_n$ be independent and identically distributed Gaussian random variables with mean $0$ and precision matrix $\uOmega\in P_{G}$. Consider a working Bayesian model that utilizes a sequence of hierarchical $G$-Wishart priors on $(\uOmega, G)$ as described in (\ref{hir_G}) and (\ref{G_prior}), satisfying Assumption \ref{asG}. Let the prior graph size cut-off value $R$ be taken such that $R\log(p_n)\to \infty$ as $n\to \infty$ and $|E_{0n}| \leq R$. Also assume that the true data-generating mechanism (see subsection \ref{subsec9.2}) satisfies Assumptions \ref{asD} and \ref{asF}. Then, if $p_n \sim n^{\delta}$ for some fixed $\delta\in (0,1)$, we have
  \begin{align*}
    \Pi_{4n}\left(\fnorm{\uOmega-\bar{\uOmega}} > M \sqrt{\frac{(p_n+ |E_{0n}|)\log p_n}{n}} \,\middle|\, \uYn\right) \overset{P}{\to} 0
  \end{align*}
  as $n \to \infty$, under $\mathbb{P}_{0,G}$ and for sufficiently large constant $M$.
\end{theorem}

\noindent The proof of Theorem \ref{postthm5} is available in the supplemental document (\cite{supp1}). Regarding the proof of Theorem \ref{postthm5}, as in \cite{Sagar,banerjee2} we rely on the general theory of posterior convergence rates developed in \cite[Theorem 2.1]{ghosal_post_conv}. Apart from this high-level idea, our argument is fundamentally different: \cite{banerjee2} and \cite{Sagar} work with i.i.d.\ Laplace and horseshoe-type priors, respectively, whereas we consider the substantially more involved $G$-Wishart prior. In this sense, our proof is completely novel, and along the way we establish several intermediate results (for example, Lemma~\href{a}{4.1} in the supplementary material) that may be of independent interest for decomposable precision matrix estimation. Moreover, to the best of our knowledge, the only related work using a $G$-Wishart prior for the precision matrix is \cite{liu2019empirical}, who study an empirical $G$-Wishart prior and obtain posterior consistency without assuming graph selection consistency. In their setting, however, only the true graph is assumed decomposable, and model misspecification is penalized via a fractional-likelihood approach, so their model setting and proof technique are entirely different from ours.

Regarding the assumption $|E_{0n}| \leq R$, recall that $R$ is the prior cut-off value (see (\ref{G_prior})) used to exclude unrealistically large graphs. This condition simply ensures that the true graph $G_0$ receives positive prior probability, which is a very mild requirement and is also imposed in establishing posterior consistency in \cite{Lee, Yabo_Niu}. The condition $R \log(p_n) \to \infty$ as $n \to \infty$ is likewise easily satisfied; for instance, it holds whenever $R = O(1)$ or $R \propto \frac{n}{\log p_n}$. The choice $R \propto \frac{n}{\log p_n}$ is standard in the high-dimensional parameter estimation literature for ruling out unrealistic graph sizes; see, for example, \cite{Lee, Lee_Lee, Ghosh_Khare}. The condition $p_n \sim n^{\delta}, \delta \in (0,1)$ has likewise been imposed in \cite{Sagar,liu2019empirical} to obtain posterior consistency without requiring minimal signal strength assumptions in the Bayesian precision matrix estimation literature.

\noindent With this important and interesting detour in place, we now return to the main line of argument for establishing the BvM-type result. With the posterior contraction results in Theorems~\ref{postthm3} and \ref{postthm4} at hand, we are now ready to present our BvM theorem for the sparse precision matrix under an unknown underlying graph $G$.
 For this purpose, we need a flatness assumption on the prior distribution for the sparse precision matrix similar to Assumption \ref{asI}. Let $\varepsilon_{4,n}= \sqrt{a^{G_{0,n}}\log p_n/n}$ be the posterior contraction rate in Theorem \ref{postthm3} and define $\rho^{\pi_{4}}(\varepsilon_{4,n})$  as
\begin{align}\label{flat}
\rho^{\pi_{4}}(\varepsilon_{4,n}):= \underset{\mathbf{T_4}\in D(\varepsilon_{4,n})}{\sup}\;\left\lvert\frac{\pi_{4}(\bar{\uOmega}+\dfrac{\uT_4}{\sqrt{n}})}{\pi_{4}(\bar{\uOmega})}-1\right\rvert,
\end{align}
where $\pi_{4}(\cdot)$ denotes the prior density for $\uOmega$, $B_{4n}=\{\uT_4:\hat{\uOmega}_{\hat{G}}+\frac{\uT_4}{\sqrt{n}}\in P_{G_0}\}$ and $D(\varepsilon_{4,n})=\{\uT_4\in B_{4n}\mid\norm{\uT_4} \leq \sqrt{n}\varepsilon_{4,n} \}$. We will now formally state the flatness assumption for the sequence of prior distributions $\{\Pi_{4n}(\cdot)\}_{n\geq1}$.
\begin{assumption}\label{asK}
    The sequence of conditional prior distributions $\{\Pi_{4n}(\cdot\mid G_0)\}_{n\geq1}$ on $\uOmega$ with support on $P_{G_0}$ is flat around true precision matrix $\bar{\uOmega}$ i.e. $\rho^{\pi_{4}}(\varepsilon_{4,n})\to 0$ in probability as $n\to\infty$.
\end{assumption}
Similar to Lemma \ref{gw}, the following lemma below shows that the $G$-Wishart prior distribution satisfies Assumption \ref{asK} under certain conditions on the hyperparameters. We will skip the proof since it is similar to the proof of Lemma \ref{gw}.
\begin{lemma}\label{gw1}
Let $W_{G_0}(\beta,\uPsi_{3n})$ denote the $G$-Wishart prior distribution defined in (\ref{prior4p}). Under Assumptions \ref{asD}, \ref{asH}, and \ref{asG}, this prior distribution satisfies Assumption \ref{asK}.
\end{lemma}

Let us now present the key result of this section. We will establish the BvM results for the precision matrix $\boldsymbol{\Omega}$ under the concentration graphical model when the true graph is unknown, as stated in the theorem below. The proof of Theorem \ref{th_omega_sp1} is included in the supplemental document (\cite{supp1}).
\begin{theorem}{\textbf{(BvM Theorem for a Sparse Precision Matrix with unknown $G$)}}\label{th_omega_sp1}
Let $Y_1, Y_2, Y_3, \cdots, Y_n$ be independent and identically distributed Gaussian random variables with mean $0$ and precision matrix $\uOmega\in P_{G}$. Consider a working Bayesian model that utilizes a hierarchical $G$-Wishart prior on $(\uOmega, G)$ as described in (\ref{hir_G}), satisfying Assumption \ref{asK} and \ref{asJ}, and assume that true data generating mechanism (see Section \ref{sec8}) satisfies Assumptions \ref{asD} and \ref{asH}. Then for any consistent estimator $\hat{G}_n$ of $G$ we have
  \begin{align*}
      \int_{B_{4n}} \lvert \pi_{4n}(\uT_4\mid\uY^n)-\phi(\uT_4;\hat{\uOmega}_{\hat{G}_n})\uone\{\uT_4\in M_{\hat{G}_n}\}\rvert\; d\uT_4 \overset{P}{\to}0,\; \textit{as}\;n\to\infty\; 
  \end{align*}
  where $\uT_4=\sqrt{n}(\uOmega-\hat{\uOmega}_{\hat{G}_n})$ and $\phi(\uT_4;\hat{\uOmega}_{\hat{G}_n})$ denotes the probability density function of the $\mathcal{SSMN}_{\hat{G}_n}(\\\uO,\;2 \uD_{\hat{G}_n}^T(\hat{\uOmega}_{\hat{G}_n}\otimes\hat{\uOmega}_{\hat{G}_n})\uD_{\hat{G}_n})$ distribution as defined in Section \ref{sec7}.
\end{theorem}

\noindent Note that Theorem \ref{th_omega_sp1} states that for large $n$, we can approximate the posterior distribution of $\sqrt{n}(\boldsymbol{\Omega} - \hat{\boldsymbol{\Omega}}_{\hat{G}_n})$ by a $\mathcal{SSMN}_{\hat{G}_n}(\mathbf{0},\;2 \mathbf{D}_{\hat{G}_n}^T(\hat{\boldsymbol{\Omega}}_{\hat{G}_n} \otimes \hat{\boldsymbol{\Omega}}_{\hat{G}_n}) \mathbf{D}_{\hat{G}_n})$ distribution, even when the true underlying graph is unknown. All we need is a consistent estimator of the true graph $\hat{G}_n$, and this estimator $\hat{G}_n$ does not need to be decomposable. Additionally, the distribution $\mathcal{SSMN}_{\hat{G}_n}(\mathbf{0},\;2 \mathbf{D}_{\hat{G}_n}^T(\hat{\boldsymbol{\Omega}}_{\hat{G}_n} \otimes \hat{\boldsymbol{\Omega}}_{\hat{G}_n}) \mathbf{D}_{\hat{G}_n})$ is entirely data-driven and thus serves the purposes of the Bernstein von-Mises theorem.
\begin{remark}
    From a purely theoretical standpoint, use of $\hat{G}_n$ is not essential: in the BvM theorem, $\hat{G}_n$ could be replaced by $G_0$, in which case no assumption on the existence or consistency of $\hat{G}_n$ is needed. Presenting BvM-type results in terms of the true underlying parameter is standard in the literature; see, for example, \cite{FangGhosh2024, castillo1, RyanNing2020}. Our data-driven formulation simply mirrors the practical setting while maintaining full theoretical validity.
\end{remark}

\section{Equivalence of different norms in terms of convergence}\label{sec11}

As previously mentioned, the use of the total variation (TV) norm is not exclusive to this problem. In this section, we will demonstrate that similar results to those in Theorem \ref{th_sigma}, \ref{th_omega}, and \ref{th_omega_sp} can be obtained by considering alternative norms.

We consider two densities, namely, $f_n$ and $g_n$, both of which are absolutely continuous with respect to a $\sigma$-finite measure $\mu$ that depends on $n$. We can define the $\alpha$-divergence, proposed by \cite{renyi}, between $f_n$ and $g_n$ as follows

\begin{align}\label{renyi}
R_{\alpha}(f_n,\;g_n)=\frac{1}{\alpha-1}\log\left[\int f_n^{\alpha}g_n^{1-\alpha}d\mu\right].
\end{align}
Similarly, we can define another type of divergence, denoted as $D_{\alpha}$ (or information divergence of type $(1-\alpha)$), given by
\begin{align}\label{d_alpha}
D_{\alpha}(f_n,g_n)=\frac{1}{\alpha(1-\alpha)}\left[1-\int f_n^{\alpha}g_n^{1-\alpha}d\mu\right].
\end{align}
It is evident that
\begin{align}\label{renyi_d_alpha}
R_{\alpha}(f_n,g_n)=\frac{1}{\alpha-1}\log\left[1-\alpha(1-\alpha)D_{\alpha}(f_n,g_n)\right].
\end{align}
Additionally, as a special case of the latter, we have $D_{1/2}(f_n,g_n)= 2H^2(f_n,g_n)$, where $H(f_n,g_n)=\left\{\int(f_n^{1/2}-g_n^{1/2})^2 d\mu\right\}^{1/2}$ represents the Bhattacharya-Hellinger distance between the densities $f_n$ and $g_n$ (\cite{bhattacharyya1946measure, Hellinger}). We now establish an inequality between the total variation distance $TV(f_n,g_n)$ and $D_{\alpha}(f_n,g_n)$. The following lemma is due to \cite{ghosh}.
\begin{lemma}\label{normequiv1}
For $0\leq\alpha\leq 1$, $\alpha(1-\alpha)D_{\alpha}(f_n,g_n)\leq\mbox{TV}(f_n,g_n)$. 
\end{lemma}

This result shows that if $TV(f_n,g_n)\rightarrow 0$, then $D_{\alpha}(f_n,g_n)$ also tends to 0 for all $\alpha\in(0,1)$. Additionally, the Hellinger divergence measure yields the inequality $H^2(f_n,g_n)\leq 2\mbox{TV}(f_n,g_n)$. 
There is another result, attributed to Le Cam and presented in \cite{wainwright} as an exercise, that provides an upper bound for $TV(f_n,g_n)$ in terms of $H(f_n,g_n)$ is given below
\begin{lemma}\label{normequiv2}
$[\mbox{TV}(f_n,g_n)]^2\leq H^2(f_n,g_n)\left[1-\frac{1}{4}H^2(f_n,g_n)\right]
\leq H^2(f_n,g_n)$.
\end{lemma}
Hence, Lemmas \ref{normequiv1} and \ref{normequiv2} have an important consequence, establishing the following convergence equivalence:
\begin{align}\label{normequiv3}
 H(f_n,g_n)\rightarrow 0\equiv TV(f_n,g_n)\rightarrow 0\equiv D_{\alpha}(f_n,g_n)
\rightarrow 0 \equiv R_{\alpha}(f_n,g_n) \rightarrow 0,
\end{align}
for all $0<\alpha<1$ as $n\to\infty$.
The equivalence between TV and Hellinger distances mentioned above is also stated in \cite{gibbs}, but it does not discuss the general R\'enyi or $\alpha$ divergence. Further details on other available norms and convergences can be found in \cite{ghosh}.

By utilizing the convergence equivalence stated in (\ref{normequiv3}), we can infer that in Theorems \ref{th_sigma}, \ref{th_omega}, and \ref{th_omega_sp}, the TV norm can be substituted with the Hellinger distance, general R\'enyi divergence, or $\alpha$ divergence for any $0<\alpha<1$. This extension allows for a wider range of norms to be applied, thereby broadening the scope of our results.

\section{Discussion}\label{sec11}

This article focuses on establishing high-dimensional Bernstein-von Mises (BvM) results for covariance and precision matrices within an independent and identically distributed Gaussian framework. In the unstructured setting, we establish BvM for $\uSigma$ (and for $\uOmega$) under mild regularity assumptions on several variables, the true data-generating mechanism, and for a general class of priors (Theorem \ref{th_sigma} and Theorem \ref{th_omega}). Next, we consider concentration graphical models where sparsity is introduced in the precision matrix to reduce the effective number of parameters. For this particular model, we initially improved the posterior contraction rates for the sparse $\boldsymbol{\Omega}$ under mild regularity assumptions on the number of variables and the true data-generating mechanism, as well as on the priors (see Theorems \ref{postthm1}, \ref{postthm2}, \ref{postthm3}, and \ref{postthm4}) for both cases when the true underlying graph is known or unknown. Additionally, we established Bernstein-von Mises (BvM) results for such models (Theorems \ref{th_omega_sp} and \ref{th_omega_sp1}).

Another common approach to introduce a low-dimensional structure in the covariance matrix is to induce sparsity in the Cholesky parameter of the precision matrix (rather than the precision matrix itself). The sparsity patterns in these matrices can be uniquely represented using appropriately directed graphs, leading to models known as directed acyclic graph models \cite{cao, geiger, smith}. However, it should be noted that without the assumption of decomposability, the precision matrix and its Cholesky parameter are not guaranteed to share the same sparsity structure. Hence, establishing BvM results for a general directed acyclic graph model and more complex covariance structures remains an open problem.  Nonetheless, Theorem \ref{th_omega_sp} and \ref{th_omega_sp1} signify a promising step forward in this direction.

\begin{acks}[Acknowledgments]
The authors are grateful to the two anonymous referees and the Associate Editor for their constructive feedback, which greatly enhanced the quality of this paper.
\end{acks}

\begin{funding}
Kshitij Khare's work for this paper was supported by NSF-DMS-2410677.
\end{funding}
\begin{supplement}
\stitle{Supplement to "High-dimensional Bernstein Von-Mises theorems for covariance and precision matrices"}
\sdescription{The supplement (\cite{supp1}) provides the remaining proofs.}
\end{supplement}


\bibliographystyle{imsart-nameyear} 
\bibliography{bibfile}       


\end{document}